\documentclass[a4paper,12pt]{article}
\usepackage{currfile,datetime}
\usepackage{cite}
\usepackage{amsthm}
\usepackage{amsmath}
\usepackage{amssymb}
\usepackage{bbm}
\usepackage{geometry}
\usepackage{epsfig}
\usepackage{listings}
\usepackage{paralist}
\usepackage{enumerate}
\usepackage{comment}
\usepackage{nicefrac}
\usepackage{verbatim}
\usepackage{multirow}
\usepackage{xcolor}

\usepackage{mathrsfs}  
\usepackage{mathtools}
\usepackage{xparse}

\usepackage[colorlinks,linkcolor=red,citecolor=green]{hyperref}

\usepackage[capitalise,compress]{cleveref} 

\crefname{equation}{}{}

\newtheorem{lemma}{Lemma}[section]

\newtheorem{theorem}[lemma]{Theorem}

\crefname{subsection}{Subsection}{Subsections}
\crefname{enumi}{item}{items}
\creflabelformat{enumi}{(#2\textup{#1}#3)}

\newcommand{\vast}[2]{\left#2 \rule{0pt}{#1}\kern-.25ex\right.}

\newcommand{\F}{\mathcal{F}}
\newcommand{\R}{\mathbbm{R}}
\newcommand{\Rd}{\mathbbm{R}^d}
\newcommand{\Rdd}{\mathbbm{R}^{d\times d}}
\newcommand{\Rm}{\mathbbm{R}^m}
\newcommand{\Rdm}{\mathbbm{R}^{d\times m}}
\newcommand{\N}{\mathbbm{N}}

\renewcommand{\P}{\mathbbm{P}}
\newcommand{\E}{\mathbbm{E}}
\newcommand{\tF}{\text{F}}

\newcommand{\funcH}[2]{{\left\vert\kern-0.25ex\left\vert\kern-0.25ex\left\vert #1     \right\vert\kern-0.25ex\right\vert\kern-0.25ex\right\vert}_{2,#2}}

\newcommand{\funcN}[2]{{\left\vert\kern-0.25ex\left\vert\kern-0.25ex\left\vert #1     \right\vert\kern-0.25ex\right\vert\kern-0.25ex\right\vert}_{1,#2}}

\newcommand{\icoi}{\ensuremath{[0,\infty)}}
\newcommand{\ioci}{\ensuremath{(0,\infty]}}
\newcommand{\iooi}{\ensuremath{(0,\infty)}}

\newcommand{\icc}[2]{\ensuremath{\mathopen[ #1,#2  \mathclose] }}

\newcommand{\Hess}{\text{Hess}\,}
\newcommand{\frace}{\ensuremath{\frac{1}{2}}}
\newcommand{\tfrace}{\ensuremath{\tfrac{1}{2}}}
\newcommand{\fracv}{\ensuremath{\frac{1}{4}}}

\newcommand{\tint}[2]{\textstyle \int_{#1}^{#2} \displaystyle}
\newcommand{\oo}[1]{\ensuremath{( #1 ) }}
\newcommand{\bigoo}[1]{\ensuremath{\big( #1 \big) }}
\newcommand{\Bigoo}[1]{\ensuremath{\Big( #1 \Big) }}
\newcommand{\biggoo}[1]{\ensuremath{\bigg( #1 \bigg) }}
\newcommand{\Biggoo}[1]{\ensuremath{\Bigg( #1 \Bigg) }}
\newcommand{\lroo}[1]{\ensuremath{\left( #1 \right) }}

\newcommand{\oc}[1]{\ensuremath{( #1 ] }}

\newcommand{\cc}[1]{\ensuremath{[ #1 ] }}
\newcommand{\bigcc}[1]{\ensuremath{\big[ #1 \big] }}
\newcommand{\Bigcc}[1]{\ensuremath{\Big[ #1 \Big] }}
\newcommand{\biggcc}[1]{\ensuremath{\bigg[ #1 \bigg] }}
\newcommand{\Biggcc}[1]{\ensuremath{\Bigg[ #1 \Bigg] }}
\newcommand{\lrcc}[1]{\ensuremath{\left[ #1 \right] }}

\newcommand{\co}[1]{\ensuremath{[ #1 ) }}

\newcommand{\bb}[1]{\ensuremath{| #1 | }}
\newcommand{\bigbb}[1]{\ensuremath{\big| #1 \big| }}
\newcommand{\Bigbb}[1]{\ensuremath{\Big| #1 \Big| }}
\newcommand{\biggbb}[1]{\ensuremath{\bigg| #1 \bigg| }}

\newcommand{\lrbb}[1]{\ensuremath{\left| #1 \right| }}

\newcommand{\gk}[1]{\ensuremath{\{ #1 \} }}
\newcommand{\biggk}[1]{\ensuremath{\big\{ #1 \big\} }}

\newcommand{\norm}[1]{\ensuremath{\| #1 \|}}
\newcommand{\bignorm}[1]{\ensuremath{\big\| #1 \hspace{0,4mm}\big\|}}
\newcommand{\Bignorm}[1]{\ensuremath{\Big\| #1 \Big\|}}
\newcommand{\biggnorm}[1]{\ensuremath{\bigg\| #1 \bigg\|}}

\newcommand{\lrnorm}[1]{\ensuremath{\left\| #1 \right\|}}
\newcommand{\normf}[1]{\ensuremath{\| #1 \|}_F}
\newcommand{\bignormf}[1]{\ensuremath{\big\| #1 \hspace{0,4mm}\big\|}_F}
\newcommand{\Bignormf}[1]{\ensuremath{\Big\| #1 \Big\|}_F}

\newcommand{\normlrd}[1]{\ensuremath{\| #1 \|}_{L(\Rd)}}
\newcommand{\bignormlrd}[1]{\ensuremath{\big\| #1 \hspace{0,4mm}\big\|}_{L(\Rd)}}
\newcommand{\Bignormlrd}[1]{\ensuremath{\Big\| #1 \Big\|}_{L(\Rd)}}

\newcommand{\bignormlrdrm}[1]{\ensuremath{\big\| #1 \hspace{0,4mm}\big\|}_{L(\Rd,\Rm)}}

\newcommand{\lr}[2]{\ensuremath{\langle #1 \,,\, #2 \rangle }}
\newcommand{\lrr}[2]{\ensuremath{\left\langle #1 \,,\, #2 \right\rangle }}
\newcommand{\biglr}[2]{\ensuremath{\big\langle #1 \,,\, #2 \big\rangle }}
\newcommand{\Biglr}[2]{\ensuremath{\Big\langle #1 \,,\, #2 \Big\rangle }}

\newcommand{\lrrf}[2]{\ensuremath{\left\langle #1 \,,\, #2 \right\rangle_F }}

\newcommand{\lrfloor}[1]{\ensuremath{\lfloor #1 \rfloor }}

\newcommand{\maxofv}[2]{ #1 \;\vee\; #2}
\newcommand{\minofv}[2]{ #1 \;\wedge\; #2}
\newcommand{\minofvttn}{\minofv{t\!}{\!\tau_N}}
\newcommand{\minofvtabtn}{\minofv{\lrfloor{t}_N\!}{\!\tau_N}}
\newcommand{\minminoft}{\minofv{\minofv{t\!}{\!\delta_n}\!}{\!\tau^N}}
\newcommand{\minminoftab}{\minofv{\minofv{\lrfloor{t}_N\!}{\!\delta_n}\!}{\!\tau^N}}

\newcommand{\tfpartial}[1]{\ensuremath{\tfrac{\partial}}{\partial #1}}

\newcommand{\LzweiOLRm}{L^2(\Omega;L(\Rm))}

\renewcommand{\-}{\!-\!}
\newcommand{\+}{\!+\!}

\NewDocumentCommand{\fabs}{sO{}m}{%
  {\IfBooleanTF{#1}
    {\fabsaux{\left|}{\right|}{#3}}
    {\fabsaux{#2|}{#2|}{#3}}}
}
\makeatletter
\newcommand{\fabsaux}[3]{\mathpalette\fabsaux@i{{#1}{#2}{#3}}}
\newcommand{\fabsaux@i}[2]{\fabsaux@ii#1#2}
\newcommand{\fabsaux@ii}[4]{%
  \sbox\z@{$\m@th#1#2#4#3$}%
  \sbox\tw@{$\m@th\|$}%
  \mathopen{\hbox to\wd\tw@{\hss\vrule height \ht\z@ depth \dp\z@ width .3\wd\tw@\hss}}%
  \mkern-2mu #4 \mkern-2mu 
  \mathclose{\hbox to\wd\tw@{\hss\vrule height \ht\z@ depth \dp\z@ width .3\wd\tw@\hss}}%
}
\makeatother

\NewDocumentCommand{\ffabs}{som}{%
  {\IfBooleanTF{#1}
    {\fabsaux{\left|}{\right|}{#3}}
    {\IfNoValueTF{#2}
      {\fnabsaux{|}{|}{#3}}
      {\fabsaux{#2|}{#2|}{#3}}
    }
  }
}
\makeatletter
\newcommand{\fnabsaux}[3]{\mathpalette\fnabsaux@i{{#1}{#2}{#3}}}
\newcommand{\fnabsaux@i}[2]{\fnabsaux@ii#1#2}
\newcommand{\fnabsaux@ii}[4]{%
  \sbox\z@{$\m@th#1#2#4#3$}%
  \sbox\tw@{$\m@th\|$}%
  \mathopen{\hbox to\wd\tw@{\hss\vrule height .8\ht\z@ depth .5\dp\z@ width .3\wd\tw@\hss}}%
  \mkern-2mu #4 \mkern-2mu 
  \mathclose{\hbox to\wd\tw@{\hss\vrule height .8\ht\z@ depth .5\dp\z@ width .3\wd\tw@\hss}}%
}
\makeatother

\title{Stopped Brownian-increment tamed Euler method}
\author{By Martin Hutzenthaler$^1$ and Kai Kisker$^2$ \\   
\small{$^1$ Faculty of Mathematics, University of Duisburg-Essen,}\\
\small{Essen, Germany; e-mail: \texttt{martin.hutzenthaler}\textcircled{\texttt{a}}\texttt{uni-due.de}}\\
\small{$^2$ Faculty of Mathematics, University of Duisburg-Essen,}\\
\small{Essen, Germany; e-mail: \texttt{kai.kisker}\textcircled{\texttt{a}}\texttt{uni-due.de}}}

\begin{document}
\newgeometry{left=18mm, right=18mm, top=20mm, bottom=30mm}
\maketitle
\makeatletter
\let\@Xmakefnmark\@makefnmark
\let\@Xthefnmark\@thefnmark
\let\@makefnmark\relax
\let\@thefnmark\relax
\@footnotetext{\emph{AMS 2010 subject classification:} 65C30}
\@footnotetext{\emph{Keywords and phrases:}
Brownian-increment tamed Euler-Maruyama method, stochastic differential equations, exponential moments, local monotonicity condition
  }
\let\@makefnmark\@Xmakefnmark
\let\@thefnmark\@Xthefnmark
\makeatother

\begin{abstract}
	In this article we propose a new explicit Euler-type approximation
	method for stochastic differential equations (SDEs).
	In this method, Brownian increments in the recursion of the Euler method
	are replaced by suitable bounded functions of the Brownian increments.
	We prove strong convergence rate one-half for a large class of SDEs
	with polynomial coefficient functions whose local monotonicity constant
	grows at most like the logarithm of a Lyapunov-type function.
\end{abstract}

\section{Introduction}\label{sec:intro}

Stochastic differential equations (SDEs) are abundant in applications.
Solutions are typically not known explicitly and need to approximated numerically.
The classical approximation method was introduced in 1955 by Maruyama; see \cite{Maruyama1955}.
Strong approximation rates have long been known in the case of globally Lipschitz continuous coefficient functions;
see, e.g., \cite{Kloeden1992}.
However, if the coefficient functions grow super-linearly, then the classical Euler-Maruyama method diverges in the strong and weak sense; see \cite{doi:10.1098/rspa.2010.0348,10.1214/12-AAP890}.

To overcome this issue, \cite{10.1214/11-AAP803} introduced so-called tamed Euler approximations
which satisfy on a uniform grid with $N\in\N$
subintervals the recursion
\begin{equation}  \begin{split}\label{eq:tamedEuler}
		\forall k\in\{0,\ldots,N-1\}\colon \quad Y_{\frac{(k+1)T}{N}}^N
		=Y_{\frac{kT}{N}}^N+\tfrac{\mu\big(Y_{\frac{kT}{N}}^N\big)}{1+\big\|\mu\big(Y_{\frac{kT}{N}}^N\big)\big\|\frac{T}{N}}\tfrac{T}N
		+\sigma\big(Y_{\frac{kT}{N}}^N\big)\Big(W_{\frac{(k+1)T}{N}}-W_{\frac{kT}{N}}\Big).
\end{split}     \end{equation}
These approximations converge with strong rate $1/2$ if the drift coefficient $\mu$ is continuously differentiable and its derivative
is bounded from above and grows at most polynomially
and
the diffusion coefficient $\sigma$ is globally Lipschitz continuous; see \cite{10.1214/11-AAP803}.
A cosmetic chink of the tamed Euler approximations in \eqref{eq:tamedEuler} is that the coefficient functions are modified, namely
$\mu$ becomes $\mu/(1+\|\mu\|T/N)$ in the $N$-th approximation.
This was improved by further tamed Euler methods as follows.
The increment-tamed Euler method in \cite[Theorem 3.15]{HutzenthalerJentzen2015} leaves the 
coefficients unchanged on with $N$ growing subsets of the state space and rescales them outside  of these
large sets.
The stopped Euler method in \cite{LiuMao2013} (generalized in \cite{HutzenthalerJentzenWang2018})
also leaves the coefficients unchanged on with $N$ growing subsets of the state space, and sets the
coefficients to zero outside of these large sets.
This latter method is appealing from a traditional point of view since it can be formulated
in terms of the classical Euler-Maruyama method.
Also the truncated Euler method of \cite{mao2015truncated} and the projected Euler method of
\cite{beyn2016stochastic} agree with the classical Euler-Maruyama method
on with $N$ growing subsets of the state space.
For further approximation methods and strong convergence rates
in the case of superlinearly growing coefficient functions
see, e.g.,
\cite{fang2016adaptive,mao2016convergence,sabanis2013note,sabanis2016Euler,wang2013tamed,zong2014convergence}.

All strong convergence results in the above mentioned literature assume
at least a global monotonicity property of the coefficients
(which implies for differentiable $\mu$ that $\mu'$ is bounded from above).
Many (or even most) interesting multi-dimensional SDEs from applications
do not have this property; e.g., stochastic Lotka-Volterra equations, stochastic Duffing-van der Pol oscillator, stochastic SIR models and the examples in \cite{cox2022local}.
Many of these SDEs have coefficients with a certain local monotonicity property
and have sufficient exponential integrability properties so that the solutions
are locally Lipschitz continuous in the initial value with respect to $L^p$;
see \cite{cox2022local,hudde2021stochastic}.
If a numerical approximation method has the same exponential integrability property
uniformly in $N$, then this method converges with rate $1/2$ to the exact solution
under mild additional assumptions; see \cite{hutzenthaler2020perturbation}.
In particular, suitably stopped and increment-tamed Euler approximations
have uniformly bounded exponential moments under suitable assumptions;
see \cite[Theorem 1.1]{HutzenthalerJentzenWang2018}
and, e.g., \cite{cozma2016exponential,jentzen2018exponential,brehier2019strong}.

An interesting question is whether one of the approximation methods
which do not modify the coefficient functions
on with $N$ growing subsets
also has uniformly bounded exponential moments under suitable assumptions.
Unfortunately, the stopped Euler method does not have this property;
see \cite[Lemma 5.1]{HutzenthalerJentzenWang2018}.
The main reason for this are the Brownian increments and the fact that
for every standard normally distributed real-valued 
random variable $V$ and every $q\in(2,\infty)$
it holds that
$\E[\exp(|V|^q)]=\infty$.
To overcome this issue, we propose in this article
a stopped Brownian-increment tamed Euler method; see \eqref{01} below.
In this method, Brownian increments are replaced by a suitable bounded function
of the Brownian increment;
see \eqref{01}.
This approach of
taming Brownian increments is new.
Moreover, we prove that this method has strong convergence rate $1/2$
for many interesting SDEs from the literature.
The following theorem illustrates the main results of this article.
\begin{theorem}\label{thm:intro}
	Let $d,m\in\N$, 
	$T,c\in\iooi$, 
	let $\mu\colon\Rd\to\Rd$, 
	$\sigma\colon\Rd\to\Rdm$
	be locally Lipschitz continuous functions with at most polynomially growing Lipschitz constants,
	let $U\in C^2(\Rd,\icoi)$,
	assume that $\Hess(U)$ grows at most polynomially,
	let $(\Omega,\F,\P,(\mathbb{F}_t)_{t\in\cc{0,T}})$ be a filtered probability space, 
	let $W\colon\icc{0}{T}\times\Omega\to\mathbbm{R}^m$ be a standard $(\mathbb{F}_t)_{t\in\cc{0,T}}$-Brownian motion with continuous sample paths,
	let $X\colon\icc{0}{T}\times\Omega\to\Rd$ be an
	adapted stochastic processes with continuous sample paths which satisfies
	that for all $t\in\icc{0}{T}$ it holds a.s.\ that
	\begin{align}
		\label{eq:X.intro}
		X_t&=X_0+\int_0^t\mu(X_s)ds+\int_0^t\sigma(X_s)dW_s,
	\end{align}
	for all $N\in\N$ let $Y^N\colon\cc{0,T}\times\Omega\to\Rd$ satisfy
	for all $k\in\gk{0,1,\dots,N-1}$ that
	$Y^N_0=X_0$ and
	\begin{align}
		\label{01}
		Y_{\frac{(k+1)T}{N}}^N&=Y^N_{\frac{kT}{N}}+\mathbbm{1}_{\norm{Y^N_{\frac{kT}{N}}}\leq e^{\sqrt{\bb{\log(T/N)}}}}
		\cdot\Bigoo{\mu\bigoo{Y^N_{\frac{kT}{N}}}\tfrac{T}{N}
			+\sigma\bigoo{Y^N_{\frac{kT}{N}}}
			\Bigcc{\tfrac{\bigoo{W_t-W_{\frac{kT}{N}}}_i}
				{\exp\bigoo{ N\bigoo{W_t-W_{\frac{kT}{N}}}_i^4/\,T}}}_{i\in\gk{1,\dots,m}}},
	\end{align}
	assume $\E\bigcc{e^{U(X_0)}}<\infty$, and assume for all $x,y\in\Rd$ that
  $\frac{1}{c}\norm{x}^{1/c}\leq1+\bb{U(x)}$ and
	\begin{align}
		\label{generator}\biglr{(\nabla U)(x)}{\mu(x)}+\tfrace\biglr{\sigma(x)}{(\Hess U)(x)\sigma(x)}_{\tF}
		+\tfrace&\norm{\oo{\sigma(x)}^*(\nabla U)(x)}^2
		\leq c\cdot U(x)\text{,}\\
		\label{05}\frac{\lr{x-y}{\mu(x)-\mu(y)}+\frac{3(1+1/c)}{2}\normf{\sigma(x)-\sigma(y)}^2}{\norm{x-y}^2}
		&\leq\!c\!+\!\frac{\bb{U(x)}+\bb{U(y)}}{8Te^{c T}}\text{.}
	\end{align}
	Then there exists a real number $C\in\R$ such that for all $N\in\N$ it holds that
	\begin{equation}\label{eq:intro.rate}
		\sup_{k\in\gk{0,1,\dots,N}}\lrnorm{X_{\frac{kT}{N}}-Y^N_{\frac{kT}{N}}}_{L^2(\Omega;\Rd)}\leq
    \tfrac{C}{\sqrt{N}}.
	\end{equation}
\end{theorem}
\noindent
Theorem \ref{thm:intro} follows directly from Theorem \ref{lemma 4}\footnote{applied for all $N\in\N$
with $p,q_0\curvearrowleft4$, $q_1\curvearrowleft\infty$, $r\curvearrowleft2$ $\overline{U}\curvearrowleft0$,
$\Pi^N\curvearrowleft\big(\lroo{x_i\,\exp\lroo{-Nx_i^4/T}}_{i\in\gk{1,\dots,m}}\big)_{x\in\R^m}$}.

In the formulation of Theorem \ref{thm:intro},
the process $X$ is a solution of the SDE \eqref{eq:X.intro} and the processes $(Y^N)_{N\in\N}$ are 
stopped Brownian-increment tamed Euler approximations of $X$. 
The taming-functions $\R\ni x\mapsto x\cdot\exp(-Nx^4/T)\in\R$, which are applied
coordinatewise to the Brownian increments, could also be replaced by other $C^2$-functions
satisfying \eqref{60}--\eqref{62}; cf.\ \eqref{2222} for a more general 
stopped Brownian-increment tamed method.
Inequality \eqref{eq:intro.rate} is the statement
that the squared $L^2$-error is asymptotically bounded by a multiple of the
(maximal) grid width.
To achieve this we assume
that the local monotony condition \eqref{05} holds
and
that the exponential moment $\E[e^{U(X_0)}]<\infty$
of the starting point is finite.
The local monotonicity constant is assumed to grow (essentially) like $U$
which satisfies the exponential Lyapunov-type condition \eqref{generator}
(which implies that $e^{U}$ is a Lyapunov-type function).
Formally the case of $U\equiv 0$ corresponds to the global monotonicity condition.   

The article is structured as follows. In Section 2 we show that
stopped Brownian-increment tamed Euler approximations are It\^o processes.
In section 3 we derive exponential moment estimates and in Section 4 we derive
moment estimates. Section 5 establishes the strong convergence rate stated in
Theorem \ref{thm:intro}.

\section{Stopped Brownian-increment tamed Euler approximations are Itô processes}
The following lemma shows that stopped Brownian-increment tamed Euler approximations
are Itô processes. For this we require that the taming functions $\Pi^N$, $N\in\N$,
are twice continuously differentiable.
\begin{lemma}[Approximations as It\^o processes]
	\label{lemma 5}
	Let $T\in\oo{0,\infty}$,
	$d,m\in\N$, 
	let $\mu\colon\Rd\to\Rd$, 
	$\sigma\colon\Rd\to\Rdm$
	be locally Lipschitz continuous functions with at most polynomially growing Lipschitz constants,
	for all $N\in\N$ let $\Pi^N\in C^2(\Rm,\Rm)$ satisfy $\Pi^N(0)=0$,
	let $(\Omega,\F,\P,(\mathbb{F}_t)_{t\in\cc{0,T}})$ be a filtered probability space,
	let $W\colon\icc{0}{T}\times\Omega\to\R^m$ be a standard $(\mathbb{F}_t)_{t\in\cc{0,T}}$-Brownian motion with continuous sample paths,
	for all $N\in\N$ let $Y^N\colon\cc{0,T}\times\Omega\to\Rd$ satisfy that $Y_0$ is $\mathbb{F}_0$-measurable
	and for all $k\in\gk{0,1,\dots,N-1}$, $t\in\bigcc{\frac{kT}{N},\frac{(k+1)T}{N}}$ that
	\begin{align}
		\label{2222}Y_t^N=Y^N_{\frac{kT}{N}}+\mathbbm{1}_{\norm{Y^N_{\frac{kT}{N}}}\leq e^{\sqrt{\bb{\log(N/T)}}}}
		\cdot\Bigoo{\mu\bigoo{Y^N_{\frac{kT}{N}}}\bigoo{t-\tfrac{kT}{N}}
			+\sigma\bigoo{Y^N_{\frac{kT}{N}}}\Pi^N\bigoo{W_t-W_{\frac{kT}{N}}}}\text{.}
	\end{align}
	Then for all $N\in\N$, $t\in\cc{0,T}$ it holds a.s. that
	\begin{equation}
		\label{24}
		\begin{aligned}
			Y_t^N=Y^N_{0}\!&+\!\int_0^t\mathbbm{1}_{s \leq\tau^N}\cdot\Bigoo{\mu(Y^N_{\lrfloor{s}_N})
				\!+\!\tfrace\sigma(Y^N_{\lrfloor{s}_N})\Delta\Pi^N\oo{W_s\!-\!W_{\lrfloor{s}_N}}}ds\\
			&+\!\int_0^t\mathbbm{1}_{s \leq\tau^N}\cdot\sigma(Y^N_{\lrfloor{s}_N})D\Pi^N\oo{W_s\!-\!W_{\lrfloor{s}_N}}dW_s
			\text{.}
		\end{aligned}
	\end{equation}
\end{lemma}
\noindent
\textit{Proof of Lemma \ref{lemma 5}}. 
Itô's formula (applied for all $N\in\N$ with $\Pi^N\in C^2(\R^m,\R^m)$
and 
$W$) and the fact that $\Pi^N(0)=0$ yield that for all $N\in\N$, $k\in\gk{0,1,\dots,N-1}$, $t\in\bigcc{\frac{kT}{N},\frac{(k+1)T}{N}}$ it holds a.s. that 
\begin{equation}
	\begin{aligned}
		Y_t^N=Y^N_{\lrfloor{t}_N}\!&+\!\int_{\frac{kT}{N}}^{\maxofv{(\minofv{t}{\tau^N})}{\frac{kT}{N}}}\!\!\mu(Y^N_{\lrfloor{s}_N})\!+\!\tfrace\sigma(Y^N_{\lrfloor{s}_N})\Delta\Pi^N\oo{W_s\!-\!W_{\lrfloor{s}_N}}ds\\
		&+\!\int_{\frac{kT}{N}}^{\maxofv{(\minofv{t}{\tau^N})}{\frac{kT}{N}}}\!\!\sigma(Y^N_{\lrfloor{s}_N})D\Pi^N\oo{W_s\!-\!W_{\lrfloor{s}_N}}dW_s
		\text{.}
	\end{aligned}
\end{equation}
This implies (\ref{24}). The proof of Lemma \ref{lemma 5} is thus completed. $ \hfill\square $

From now on we focus on the specific taming functions \eqref{eq:taming.functions}.
The next lemma establishes upper bounds for derivatives of this taming function.
\begin{lemma}[Taming functions]\label{l:taming.function}
	\label{lemma 1}
	Let $h \in \oo{0,\infty}$, 
	$m\in\N$,
	let $(\Omega,\F,\P)$ be a probability space,
	let $W=\oo{W^1,\dots,W^m}\colon$ $\icc{0}{h}\times\Omega\to\R^m$ be a standard Brownian motion with continuous sample paths,
	and let $\Pi\in C^2(\Rm,\Rm)$ be the function which satisfies for all $x=\oo{x_1,\dots,x_m}\in\Rm$ that
	\begin{equation}\label{eq:taming.functions}
		\Pi(x)=\lroo{x_i\,\exp\lroo{-\frac{x_i^4}{h}}}_{i\in\gk{1,\dots,m}}\text{.}
	\end{equation}
	Then it holds for all $t\in\cc{0,h}$ that
	\begin{align}
		\label{60}		&\norm{\Pi(W_t)}_{L^\infty(\Omega;\Rm)}&&\leq h^{\frac{1}{4}}\sqrt{m}\text{,}\\
		\label{61}		&\bignorm{D\Pi(W_t)-I_{\Rm}}_{L^2(\Omega;L(\Rm))}
		\!\!\!\!\!\!\!\!\!\!\!\!\!\!\!\!\!\!\!\!\!\!\!\!\!\!\!\!\!\!\!\!\!\!\!\!
		&&\leq 52h\sqrt{m}\text{,}\\
		\label{62}		&\bignorm{\Delta\Pi(W_t)}_{L^2(\Omega;\Rm)}		
		&&\leq32\sqrt{hm}\text{.}
	\end{align}
\end{lemma}
\noindent
\textit{Proof of Lemma \ref{lemma 1}}. 
The fact that $\sup_{x\in\icoi}xe^{-x^4}\leq 1$ implies for all $t\in\cc{0,h}$ that
\begin{align}
	\norm{\Pi(W_t)}_{L^\infty(\Omega;\Rm)}\leq m^{\frace}h^{\frac{1}{4}}\text{.}
\end{align}
This proves (\ref{60}). The fact that $\forall v\in\Rm\colon\norm{\text{Diag}(v)}_{L(\Rm)}=\norm{v}$, the fact that $\forall x\in\icoi\colon\bb{1-e^{-x}}\leq x$ yield for all $t\in\cc{0,h}$ that 
\begin{equation}
	\begin{aligned}
		&\norm{D\Pi(W_t)-I_{\Rm}}_{L^2(\Omega;L(\Rm))}\\	
		=\;&\biggnorm{\lroo{\exp\lroo{-\frac{\oo{W_t^i}^4}{h}}\lroo{1-\frac{4\oo{W_t^i}^4}{h}}-1}_{i\in\gk{1,\dots,m}}}_{L^2(\Omega;\Rm)}\\
		\leq\;&\biggnorm{\lroo{\lrbb{1-\exp\lroo{-\frac{\oo{W_t^i}^4}{h}}}+\lrbb{\exp\lroo{-\frac{\oo{W_t^i}^4}{h}}\frac{4\oo{W_t^i}^4}{h}}}_{i\in\gk{1,\dots,m}}}_{L^2(\Omega;\Rm)}\\
		\leq\;&\biggnorm{\lroo{\frac{\oo{W_t^i}^4}{h}+\frac{4\oo{W_t^i}^4}{h}}_{i\in\gk{1,\dots,m}}}_{L^2(\Omega;L(\Rm))}\\
		=\;&5h\biggnorm{\lroo{\frac{\oo{W_h^i}^4}{h^2}}_{i\in\gk{1,\dots,m}}}_{L^2(\Omega;\Rm)}\\
		=\;&5h\biggnorm{\biggcc{\sum_{i=1}^m\lroo{\frac{\oo{W_h^i}^4}{h^2}}^2}^{\frace}}_{L^2(\Omega;\R)}\\
		=\;&5h\sqrt{m}\biggnorm{\biggoo{\frac{W^1_h}{\sqrt{h}}}^{8}}_{L^1(\Omega;\R)}^{\frace}\\
		=\;&5h\sqrt{m}\cdot\sqrt{105}\\
		\leq\;&52h\sqrt{m}\text{.}
	\end{aligned}
\end{equation}
This proves (\ref{61}).
The triangle inequality and the fact that $\forall x\in\icoi\colon \bb{e^{-x}(x-1)}\leq1$ show for all $t\in\cc{0,h}$ that
\begin{equation}
	\begin{aligned}
		&\hspace{6,5mm}\norm{\Delta\Pi(W_t)}_{L^2(\Omega;\Rm)}\\
		&\;=\biggnorm{\biggoo{\exp\biggoo{-\frac{(W_t^i)^4}{h}}
				\biggoo{\frac{16(W_t^i)^7}{h^2}-\frac{20(W_t^i)^3}{h}}}_{i\in\gk{1,\dots,m}}}_{L^2(\Omega;\Rm)}\\
		&\;=\biggnorm{\biggoo{\exp\biggoo{-\frac{(W_t^i)^4}{h}}\biggoo{\frac{(W_t^i)^4}{h}-1}16\frac{(W_t^i)^3}{h}-4\frac{(W_t^i)^3}{h}}_{i\in\gk{1,\dots,m}}}_{L^2(\Omega;\Rm)}\\
		&\;\leq20\sqrt{h}\biggnorm{\biggoo{\frac{(W_h^i)^3}{h^{\frac{3}{2}}}}_{i\in\gk{1,\dots,m}}}_{L^2(\Omega;\Rm)}\\
		&\;=20\sqrt{hm}\biggnorm{\biggoo{\frac{W^1_h}{\sqrt{h}}}^3}^{\frace}_{L^1(\Omega;\R)}\\
		&\;=20\sqrt{hm}\cdot2\sqrt{\tfrac{2}{\pi}}\\
		&\;\leq32\sqrt{hm}\text{.}
	\end{aligned}
\end{equation}
This proves (\ref{62}). The proof of Lemma \ref{lemma 1} is thus completed. $ \hfill\square $

The following lemma estimates the temporal regularity of our approximations.
\begin{lemma}[Temporal regularity]
	\label{lemma 6}
	Let $T\in\iooi$,
	$\rho\in \icoi$, 
	$d,m\in\N$, 
	let $\mu\colon\Rd\to\Rd$, 
	$\sigma\colon\Rd\to\Rdm$
	be locally Lipschitz continuous functions with at most polynomially growing Lipschitz constants,
	let $U\in C^2(\Rd,\R)$, 
	let $\overline{U}\colon\Rd\to\R$ be a measurable function, 
	assume that $\Hess (U)$, $\overline{U}$ are at most polynomially growing functions,
	for all $N\in\N$ let $\Pi^N\in C^2(\Rm,\Rm)$ be the function which satisfies for all $x=\oo{x_1,\dots,x_m}\in\Rm$ that
	\begin{equation}
		\Pi^N(x)=\lroo{x_i\,\exp\lroo{-\frac{Nx_i^4}{T}}}_{i\in\gk{1,\dots,m}}\text{,}
	\end{equation}
	let $(\Omega,\F,\P,(\mathbb{F}_t)_{t\in\cc{0,T}})$ be a filtered probability space,
	let $W\colon\icc{0}{T}\times\Omega\to\R^m$ be a standard $(\mathbb{F}_t)_{t\in\cc{0,T}}$-Brownian motion with continuous sample paths,
	and for all $N\in\N$ let $Y^N\colon\cc{0,T}\times\Omega\to\Rd$ satisfy that $Y_0^N$ is $\mathbb{F}_0$-measurable
	and for all $k\in\gk{0,1,\dots,N-1}$, $t\in\bigcc{\frac{kT}{N},\frac{(k+1)T}{N}}$ that
	\begin{align}
		\label{02}Y_t^N=Y^N_{\frac{kT}{N}}+\mathbbm{1}_{\norm{Y^N_{\frac{kT}{N}}}\leq e^{\sqrt{\bb{\log(N/T)}}}}
		\cdot\Bigoo{\mu\bigoo{Y^N_{\frac{kT}{N}}}\bigoo{t-\tfrac{kT}{N}}
			+\sigma\bigoo{Y^N_{\frac{kT}{N}}}\Pi^N\bigoo{W_t-W_{\frac{kT}{N}}}}\text{.}
	\end{align}
	Then it holds for all $N\in\N$, $t\in\cc{0,T}$ it holds a.s. that 
	\begin{enumerate}[I]
		\item 	$\label{021}\bignorm{Y_t^N-Y^N_{\lrfloor{t}_N}}\leq2c^{p+1}\bigoo{\tfrac{T}{N}}^{\frac{7}{32}}\bigoo{T^{\frac{3}{4}}+\sqrt{m}}$
		\item 	$\label{022}\mathbbm{1}_{\norm{Y^N_{\lrfloor{t}_N}}\leq e^{\sqrt{\bb{\log(N/T)}}}}\cdot\norm{Y_t^N}^p\leq3^pc^{p}\bigoo{\tfrac{N}{T}}^{\frac{1}{32}}$
		\item	$\label{023}\mathbbm{1}_{\norm{Y^N_{\lrfloor{t}_N}}\leq e^{\sqrt{\bb{\log(N/T)}}}}
		\cdot\Bigoo{\bignorm{\mu\oo{Y_t^N}\-\mu(Y^N_{\lrfloor{t}_N})}
			+\tfrace\bignorm{\sigma(Y^N_{\lrfloor{t}_N})\Delta\Pi^N(W_t-W_{\lrfloor{t}_N})}}$\\
		$\leq\;2c^{p+1}3^{p}\bigoo{\tfrac{N}{T}}^{\frac{1}{32}}	\Bigoo{2c^{p+1}\bigoo{\tfrac{T}{N}}^{\frac{7}{32}}
			\bigoo{T^{\frac{3}{4}}+\sqrt{m}}+\tfrace\bignorm{\Delta\Pi^N(W_t-W_{\lrfloor{t}_N})}_{L(\Rm)}}$
		\item 	$\label{024}\mathbbm{1}_{\norm{Y^N_{\lrfloor{t}_N}}\leq e^{\sqrt{\bb{\log(N/T)}}}}
		\cdot\bignorm{\nabla U\oo{Y_t^N}}\Bigoo{\bignorm{\mu\oo{Y_t^N}\-\mu(Y^N_{\lrfloor{t}_N})}
			+\tfrace\bignorm{\sigma(Y^N_{\lrfloor{t}_N})\Delta\Pi^N(W_t-W_{\lrfloor{t}_N})}}$\\
		$\leq\;4c^{2p+2}3^{2p}\bigoo{\tfrac{N}{T}}^{\frac{1}{16}}	\Bigoo{2c^{p+1}\bigoo{\tfrac{T}{N}}^{\frac{7}{32}}
			\bigoo{T^{\frac{3}{4}}+\sqrt{m}}+\tfrace\bignorm{\Delta\Pi^N(W_t-W_{\lrfloor{t}_N})}_{L(\Rm)}}$
		\item 	$\label{025}\mathbbm{1}_{\norm{Y^N_{\lrfloor{t}_N}}\leq e^{\sqrt{\bb{\log(N/T)}}}}
		\cdot\Bigoo{\Bignormf{\sigma(Y^N_{\lrfloor{t}_N})\cdot\bigoo{I_{\Rm}-D\Pi^N(W_t)}}
			+\bignormf{\sigma(Y^N_s)-\sigma(Y^N_{\lrfloor{s}_N})}}$\\
		$\leq\;\bignorm{D\Pi^N(W_t)-I_{\Rm}}_{L(\Rm)}\cdot 2c^{p+1}\bigoo{\tfrac{N}{T}}^{\frac{1}{32}}
		+2c^{2p+2}3^{p}\bigoo{\tfrac{T}{N}}^{\frac{3}{16}}\bigoo{T^{\frac{3}{4}}+\sqrt{m}}$
	\end{enumerate}
\end{lemma}
\noindent
\textit{Proof of Lemma \ref{lemma 6}.} 
For all $N\in\N$ let $\lrfloor{\cdot}_{N}\colon\cc{0,T}\to\cc{0,T}$ be the function which satisfies for all $t\in\oc{0,T}$ that $\lrfloor{t}_{N}=\sup(\gk{0,\frac{T}{N},\frac{2T}{N},\dots,T}\cap\co{0,t})$ and $\lrfloor{0}_{N}=0$.
Choose $p\in\N$, $c\in\co{T^{\frac{1}{32}},\infty}$, and $N_0\in\N$ such that for all $x,y\in\Rd$, $N\in\N$ with $N\geq N_0$ it holds that 
\begin{align}
	\label{011}\norm{\mu(x)-\mu(y)}+\normf{\sigma(x)-\sigma(y)}\leq c\oo{1+\norm{x}^p+\norm{y}^p}\norm{x-y}\text{,}\\
	\label{012}\bb{\overline{U}(x)}+\normlrd{\Hess U(x)}+\norm{\nabla U(x)}+\bb{U(x)}+\norm{\mu(x)}+\normf{\sigma(x)}\leq c\oo{1+\norm{x}^p}\text{,}\\
	\label{013}e^{\sqrt{\bb{\log(N/T)}}}\leq c\bigoo{\tfrac{N}{T}}^{\frac{1}{32p}} \text{,}\\
	\label{014} c^p\bigoo{\tfrac{T}{N}}^{\frac{7}{32}}\bigoo{T^{\frac{3}{4}}+\sqrt{m}}\leq\bigoo{\tfrac{N}{T}}^{\frac{1}{32p}}\text{.}
\end{align}
The triangle inequality, (\ref{02}), (\ref{012}), Lemma \ref{lemma 1} (applied for every $t\in\cc{0,T}$ with $h\curvearrowleft T/N$, $W\curvearrowleft \oo{W_{\lrfloor{t}_N+s}-W_{\lrfloor{t}_N}}_{s\in\cc{0,T/N}}$), and (\ref{013}) imply that for all $t\in\cc{0,T}$, $N\in\N$ with $N\geq N_0$ it holds a.s. that
\begin{equation}
	\label{015}
	\begin{aligned}
		&\bignorm{Y_t^N-Y^N_{\lrfloor{t}_N}}\\
		\leq\;&\mathbbm{1}_{\norm{Y^N_{\lrfloor{t}_N}}\leq e^{\sqrt{\bb{\log(N/T)}}}}
		\cdot\Bigoo{\bignorm{\mu\bigoo{Y^N_{\lrfloor{t}_N}}}\bigoo{t\!-\!\lrfloor{t}}
			\!+\!\bignormlrdrm{\sigma\bigoo{Y^N_{\lrfloor{t}_N}}}
			\bignorm{\Pi^N\bigoo{W_t-W_{\lrfloor{t}_N}}}}\\
		\leq\;&\mathbbm{1}_{\norm{Y^N_{\lrfloor{t}_N}}\leq e^{\sqrt{\bb{\log(N/T)}}}}
		\cdot\Bigoo{c\bigoo{1+\bignorm{Y^N_{\lrfloor{t}_N}}^p}\tfrac{T}{N}
			+c\bigoo{1+\bignorm{Y^N_{\lrfloor{t}_N}}^p}\sqrt{m}\bigoo{\tfrac{T}{N}}^{\fracv}}\\
		=\;& \mathbbm{1}_{\norm{Y^N_{\lrfloor{t}_N}}\leq e^{\sqrt{\bb{\log(N/T)}}}}
		\cdot c\bigoo{1+\bignorm{Y^N_{\lrfloor{t}_N}}^p}
		\Bigoo{\bigoo{\tfrac{T}{N}}^{\frac{3}{4}}+\sqrt{m}}\bigoo{\tfrac{T}{N}}^{\fracv}\\
		\leq\;& c\bigoo{1+e^{p\sqrt{\bb{\log(N/T)}}}}
		\Bigoo{\bigoo{\tfrac{T}{N}}^{\frac{3}{4}}+\sqrt{m}}\bigoo{\tfrac{T}{N}}^{\fracv}\\
		\leq\;&c\Bigoo{c^p\bigoo{\tfrac{N}{T}}^{\frac{1}{32}}+c^p\bigoo{\tfrac{N}{T}}^{\frac{1}{32}}}
		\Bigoo{\bigoo{\tfrac{T}{N}}^{\frac{3}{4}}+\sqrt{m}}\bigoo{\tfrac{T}{N}}^{\fracv}\\
		\leq\;&2c^{p+1}\bigoo{\tfrac{T}{N}}^{\frac{7}{32}}\bigoo{T^{\frac{3}{4}}+\sqrt{m}}\text{.}
	\end{aligned}
\end{equation}
This proves \ref{021}.
The triangle inequality, (\ref{013}), (\ref{014}), and (\ref{015}) show that for all $t\in\cc{0,T}$, $N\in\N$ with $N\geq N_0$ it holds a.s. that 
\begin{equation}
	\begin{aligned}
		\label{019}
		\mathbbm{1}_{\norm{Y^N_{\lrfloor{t}_N}}\leq e^{\sqrt{\bb{\log(N/T)}}}}\cdot\norm{Y_t^N}^p
		&\leq\mathbbm{1}_{\norm{Y^N_{\lrfloor{t}_N}}\leq e^{\sqrt{\bb{\log(N/T)}}}}
		\cdot\Bigoo{\bignorm{Y_t^N-Y^N_{\lrfloor{t}_N}}+\bignorm{Y^N_{\lrfloor{t}_N}}}^p\\
		&\leq\Bigoo{2c^{p+1}\bigoo{\tfrac{T}{N}}^{\frac{7}{32}}\bigoo{T^{\frac{3}{4}}+\sqrt{m}}+e^{\sqrt{\bb{\log(N/T)}}}}^p\\
		&\leq\Bigoo{2c^{p+1}\bigoo{\tfrac{T}{N}}^{\frac{7}{32}}\bigoo{T^{\frac{3}{4}}+\sqrt{m}}+c\bigoo{\tfrac{N}{T}}^{\frac{1}{32p}}}^p\\
		&\leq\Bigoo{3c\bigoo{\tfrac{N}{T}}^{\frac{1}{32p}}}^p\\
		&= 3^pc^{p}\bigoo{\tfrac{N}{T}}^{\frac{1}{32}}\text{.}
	\end{aligned}
\end{equation}
This proves \ref{022}.
Next (\ref{011}), (\ref{012}), (\ref{013}), (\ref{015}), and (\ref{019}) imply that for all $t\in\cc{0,T}$, $N\in\N$ with $N\geq N_0$ it holds a.s. that 
\begin{equation}
	\begin{aligned}
		\label{018}
		&\mathbbm{1}_{\norm{Y^N_{\lrfloor{t}_N}}\leq e^{\sqrt{\bb{\log(N/T)}}}}
		\cdot\Bigoo{\bignorm{\mu\oo{Y_t^N}\-\mu(Y^N_{\lrfloor{t}_N})}
			+\tfrace\bignorm{\sigma(Y^N_{\lrfloor{t}_N})\Delta\Pi^N(W_t-W_{\lrfloor{t}_N})}}\\
		\leq\;&\mathbbm{1}_{\norm{Y^N_{\lrfloor{t}_N}}\leq e^{\sqrt{\bb{\log(N/T)}}}}
		\cdot\Bigoo{\bignorm{\mu\oo{Y_t^N}\-\mu(Y^N_{\lrfloor{t}_N})}
			\+\tfrace\bignorm{\Delta\Pi^N(W_t\-W_{\lrfloor{t}_N})}_{L(\Rm)}\bignormf{\sigma(Y^N_{\lrfloor{t}_N})}}\\
		\leq\;&c\mathbbm{1}_{\norm{Y^N_{\lrfloor{t}_N}}\leq e^{\sqrt{\bb{\log(N/T)}}}}
		\cdot\Bigoo{c\Bigcc{1+\bignorm{Y_t^N}^p+\bignorm{Y^N_{\lrfloor{t}_N}}^p}\bignorm{Y_t^N-Y^N_{\lrfloor{t}_N}}\\
			&\;+\tfrace\bignorm{\Delta\Pi^N(W_t-W_{\lrfloor{t}_N})}_{L(\Rm)}c\bigoo{1+\bignorm{Y_{\lrfloor{t}_N}^N}^p}}\\
		\leq\;&c\mathbbm{1}_{\norm{Y^N_{\lrfloor{t}_N}}\leq e^{\sqrt{\bb{\log(N/T)}}}}
		\cdot\Bigoo{1\!+\!\bignorm{Y_t^N}^p	\!+\!\bignorm{Y^N_{\lrfloor{t}_N}}^p}
		\Bigoo{\bignorm{Y_t^N\!-\!Y^N_{\lrfloor{t}_N}}\!+\!\tfrace\bignorm{\Delta\Pi^N(W_t-W_{\lrfloor{t}_N})}_{L(\Rm)}}\\
		\leq\;&c\!\Bigoo{\!3^pc^{p}\bigoo{\tfrac{N}{T}}^{\frac{1}{32}}	\!\!+\!2e^{p\sqrt{\bb{\log(N/T)}}}}\!
		\Bigoo{\!2c^{p+1}\!\bigoo{\tfrac{T}{N}}^{\frac{7}{32}}\bigoo{T^{\frac{3}{4}}\!+\!\sqrt{m}}\!+\!\tfrace\bignorm{\Delta\Pi^N(W_t\!-\!W_{\lrfloor{t}_N})}_{L(\Rm)}}\\
		\leq\;&c\Bigoo{3^pc^{p}\bigoo{\tfrac{N}{T}}^{\frac{1}{32}}	\!+\!2c^p\bigoo{\tfrac{N}{T}}^{\frac{1}{32}}}
		\Bigoo{2c^{p+1}\bigoo{\tfrac{T}{N}}^{\frac{7}{32}}\bigoo{T^{\frac{3}{4}}\!+\!\sqrt{m}}\!+\!\tfrace\bignorm{\Delta\Pi^N(W_t-W_{\lrfloor{t}_N})}_{L(\Rm)}}\\
		\leq\;&c\Bigoo{2c^{p}3^p\bigoo{\tfrac{N}{T}}^{\frac{1}{32}}	}\Bigoo{2c^{p+1}\bigoo{\tfrac{T}{N}}^{\frac{7}{32}}
			\bigoo{T^{\frac{3}{4}}+\sqrt{m}}+\tfrace\bignorm{\Delta\Pi^N(W_t-W_{\lrfloor{t}_N})}_{L(\Rm)}}\\
		\leq\;&2c^{p+1}3^{p}\bigoo{\tfrac{N}{T}}^{\frac{1}{32}}	\Bigoo{2c^{p+1}\bigoo{\tfrac{T}{N}}^{\frac{7}{32}}
			\bigoo{T^{\frac{3}{4}}+\sqrt{m}}+\tfrace\bignorm{\Delta\Pi^N(W_t-W_{\lrfloor{t}_N})}_{L(\Rm)}}\text{.}	
	\end{aligned}
\end{equation}
This proves \ref{023}.
This and (\ref{012}) show that for all $t\in\cc{0,T}$, $N\in\N$ with $N\geq N_0$ it holds a.s. that 
\begin{equation}
	\begin{aligned}
		\label{020}
		&\mathbbm{1}_{\norm{Y^N_{\lrfloor{t}_N}}\leq e^{\sqrt{\bb{\log(N/T)}}}}
		\cdot\bignorm{\nabla U\oo{Y_t^N}}
		\Bigoo{\bignorm{\mu\oo{Y_t^N}\-\mu(Y^N_{\lrfloor{t}_N})}
			+\tfrace\bignorm{\sigma(Y^N_{\lrfloor{t}_N})\Delta\Pi^N(W_t-W_{\lrfloor{t}_N})}}\\
		\leq\;&\mathbbm{1}_{\norm{Y^N_{\lrfloor{t}_N}}\leq e^{\sqrt{\bb{\log(N/T)}}}}
		\cdot\bigoo{1+\norm{Y_t^N}^p}\\
		&\;\cdot2c^{p+2}3^{p}\bigoo{\tfrac{N}{T}}^{\frac{1}{32}}	\Bigoo{2c^{p+1}\bigoo{\tfrac{T}{N}}^{\frac{7}{32}}
			\bigoo{T^{\frac{3}{4}}+\sqrt{m}}+\tfrace\bignorm{\Delta\Pi^N(W_t-W_{\lrfloor{t}_N})}_{L(\Rm)}}\\
		\leq\;&\mathbbm{1}_{\norm{Y^N_{\lrfloor{t}_N}}\leq e^{\sqrt{\bb{\log(N/T)}}}}
		\cdot\Bigoo{1\!+\!\bignorm{Y_t^N}^p	\!+\!\bignorm{Y^N_{\lrfloor{t}_N}}^p}\\
		&\;\cdot2c^{p+2}3^{p}\bigoo{\tfrac{N}{T}}^{\frac{1}{32}}	\Bigoo{2c^{p+1}\bigoo{\tfrac{T}{N}}^{\frac{7}{32}}
			\bigoo{T^{\frac{3}{4}}+\sqrt{m}}+\tfrace\bignorm{\Delta\Pi^N(W_t-W_{\lrfloor{t}_N})}_{L(\Rm)}}\\
		\leq\;&\Bigoo{\!3^pc^{p}\bigoo{\tfrac{N}{T}}^{\frac{1}{32}}	\!\!+\!2e^{p\sqrt{\bb{\log(N/T)}}}}\\
		&\;\cdot2c^{p+2}3^{p}\bigoo{\tfrac{N}{T}}^{\frac{1}{32}}	\Bigoo{2c^{p+1}\bigoo{\tfrac{T}{N}}^{\frac{7}{32}}
			\bigoo{T^{\frac{3}{4}}+\sqrt{m}}+\tfrace\bignorm{\Delta\Pi^N(W_t-W_{\lrfloor{t}_N})}_{L(\Rm)}}\\
		\leq\;&\Bigoo{3^pc^{p}\bigoo{\tfrac{N}{T}}^{\frac{1}{32}}	\!+\!2c^p\bigoo{\tfrac{N}{T}}^{\frac{1}{32}}}\\
		&\;\cdot2c^{p+2}3^{p}\bigoo{\tfrac{N}{T}}^{\frac{1}{32}}	\Bigoo{2c^{p+1}\bigoo{\tfrac{T}{N}}^{\frac{7}{32}}
			\bigoo{T^{\frac{3}{4}}+\sqrt{m}}+\tfrace\bignorm{\Delta\Pi^N(W_t-W_{\lrfloor{t}_N})}_{L(\Rm)}}\\
		\leq\;&\Bigoo{2c^{p}3^p\bigoo{\tfrac{N}{T}}^{\frac{1}{32}}	}
		\cdot2c^{p+2}3^{p}\bigoo{\tfrac{N}{T}}^{\frac{1}{32}}	\Bigoo{2c^{p+1}\bigoo{\tfrac{T}{N}}^{\frac{7}{32}}
			\bigoo{T^{\frac{3}{4}}+\sqrt{m}}+\tfrace\bignorm{\Delta\Pi^N(W_t-W_{\lrfloor{t}_N})}_{L(\Rm)}}\\
		\leq\;&4c^{2p+2}3^{2p}\bigoo{\tfrac{N}{T}}^{\frac{1}{16}}	\Bigoo{2c^{p+1}\bigoo{\tfrac{T}{N}}^{\frac{7}{32}}
			\bigoo{T^{\frac{3}{4}}+\sqrt{m}}+\tfrace\bignorm{\Delta\Pi^N(W_t-W_{\lrfloor{t}_N})}_{L(\Rm)}}\text{.}	
	\end{aligned}
\end{equation}
This proves \ref{024}.
Moreover (\ref{013}), (\ref{015}), and (\ref{019}) show that for all $t\in\cc{0,T}$, $N\in\N$ with $N\geq N_0$ it holds a.s. that
\begin{equation}
	\begin{aligned}
		\label{070}
		&\mathbbm{1}_{\norm{Y^N_{\lrfloor{t}_N}}\leq e^{\sqrt{\bb{\log(N/T)}}}}
		\cdot\biggcc{\Bignormf{\sigma(Y^N_{\lrfloor{t}_N})\cdot\bigoo{I_{\Rm}-D\Pi^N(W_t)}}
			+\bignormf{\sigma(Y^N_s)-\sigma(Y^N_{\lrfloor{s}_N})}}\\
		\leq\;&\mathbbm{1}_{\norm{Y^N_{\lrfloor{t}_N}}\leq e^{\sqrt{\bb{\log(N/T)}}}}
		\cdot\biggcc{\bignorm{D\Pi^N(W_t)-I_{\Rm}}_{L(\Rm)}\cdot\bignormf{\sigma(Y^N_{\lrfloor{t}_N})}\\
			&+c\Bigoo{1+\bignorm{Y_t^N}^p
				+\bignorm{Y^N_{\lrfloor{t}_N}}^p}\bignorm{Y_t^N-Y^N_{\lrfloor{t}_N}}}\\
		\leq\;&\mathbbm{1}_{\norm{Y^N_{\lrfloor{t}_N}}\leq e^{\sqrt{\bb{\log(N/T)}}}}
		\cdot\bignorm{D\Pi^N(W_t)-I_{\Rm}}_{L(\Rm)}\cdot c\bigoo{1+\bignorm{Y_{\lrfloor{t}_N}^N}^p}\\
		&+c\Bigoo{3^pc^{p}\bigoo{\tfrac{N}{T}}^{\frac{1}{32}}+2e^{p\sqrt{\bb{\log(N/T)}}}}
		2c^{p+1}\bigoo{\tfrac{T}{N}}^{\frac{7}{32}}\bigoo{T^{\frac{3}{4}}+\sqrt{m}}\\
		\leq\;&\bignorm{D\Pi^N(W_t)-I_{\Rm}}_{L(\Rm)}\cdot 2ce^{p\sqrt{\bb{\log(N/T)}}}\\
		&+c\Bigoo{3^pc^{p}\bigoo{\tfrac{N}{T}}^{\frac{1}{32}}+2c^p\bigoo{\tfrac{N}{T}}^{\frac{1}{32}}}
		2c^{p+1}\bigoo{\tfrac{T}{N}}^{\frac{7}{32}}\bigoo{T^{\frac{3}{4}}+\sqrt{m}}\\
		\leq\;&\bignorm{D\Pi^N(W_t)-I_{\Rm}}_{L(\Rm)}\cdot 2c^{p+1}\bigoo{\tfrac{N}{T}}^{\frac{1}{32}}
		+2c^{2p+2}3^{p}\bigoo{\tfrac{N}{T}}^{\frac{1}{32}}
		\bigoo{\tfrac{T}{N}}^{\frac{7}{32}}\bigoo{T^{\frac{3}{4}}+\sqrt{m}}\\
		\leq\;&\bignorm{D\Pi^N(W_t)-I_{\Rm}}_{L(\Rm)}\cdot 2c^{p+1}\bigoo{\tfrac{N}{T}}^{\frac{1}{32}}
		+2c^{2p+2}3^{p}\bigoo{\tfrac{T}{N}}^{\frac{3}{16}}\bigoo{T^{\frac{3}{4}}+\sqrt{m}}\text{.}
	\end{aligned}
\end{equation}
This proves \ref{025}. The proof of Lemma \ref{lemma 6} is thus completed. $ \hfill\square $\\

\section{Exponential moment estimates}
In this subsection we prove that
 that stopped Brownian-increment tamed Euler approximations
have uniformly bounded exponential moments. This is the core of our analysis.
Lemma \ref{lemma 2} is motivated by \cite{HutzenthalerJentzenWang2018} which proves the
respective result for stopped increment-tamed Euler approximations.
\begin{lemma}[Exponential moments]
	\label{lemma 2}
	Let $T\in\iooi$,
	$\rho\in \icoi$, 
	$d,m\in\N$, 
	let $\mu\colon\Rd\to\Rd$, 
	$\sigma\colon\Rd\to\Rdm$
	be locally Lipschitz continuous functions with at most polynomially growing Lipschitz constants,
	let $U\in C^2(\Rd,\R)$, 
	let $\overline{U}\colon\Rd\to\R$ be measurable function, 
	assume that $\Hess (U)$, $\overline{U}$ are at most polynomially growing functions,
	for all $N\in\N$ let $\Pi^N\in C^2(\Rm,\Rm)$ be the function which satisfies for all $x=\oo{x_1,\dots,x_m}\in\Rm$ that
	\begin{equation}
		\Pi^N(x)=\lroo{x_i\,\exp\lroo{-\frac{Nx_i^4}{T}}}_{i\in\gk{1,\dots,m}}\text{,}
	\end{equation}
	let $(\Omega,\F,\P,(\mathbb{F}_t)_{t\in\cc{0,T}})$ be a filtered probability space,
	let $W\colon\icc{0}{T}\times\Omega\to\R^m$ be a standard $(\mathbb{F}_t)_{t\in\cc{0,T}}$-Brownian motion with continuous sample paths,
	for all $N\in\N$ let $Y^N\colon\cc{0,T}\times\Omega\to\Rd$ satisfy that $Y_0$ is $\mathbb{F}_0$-measurable
	and for all $k\in\gk{0,1,\dots,N-1}$, $t\in\bigcc{\frac{kT}{N},\frac{(k+1)T}{N}}$ that
	\begin{align}
		\label{2}Y_t^N=Y^N_{\frac{kT}{N}}+\mathbbm{1}_{\norm{Y^N_{\frac{kT}{N}}}\leq e^{\sqrt{\bb{\log(N/T)}}}}
		\cdot\Bigoo{\mu\bigoo{Y^N_{\frac{kT}{N}}}\bigoo{t-\tfrac{kT}{N}}
			+\sigma\bigoo{Y^N_{\frac{kT}{N}}}\Pi^N\bigoo{W_t-W_{\frac{kT}{N}}}}\text{,}
	\end{align}
	for all $N\in\N$ let $\tau^N\colon\Omega\to\cc{0,T}$ be the function which satisfies that
	$\tau^N=\inf\bigoo{\biggk{\frac{kT}{N}\in\cc{0,T}\colon k\in\gk{0,1,\dots,N},\, \norm{Y_{\frac{kT}{N}}^N}>e^{\sqrt{\bb{\log(N/T)}}}}\cup\gk{T}}$,

	and assume for all $y\in\Rd$ that
	\begin{align}
		\label{1}\!\!\biglr{(\nabla U)(y)}{\mu(y)}+\tfrace\biglr{\sigma(y)}{(\Hess U)(y)\sigma(y)}_{\tF}
		+&\tfrace\norm{\oo{\sigma(y)}^*(\nabla U)(y)}^2\!+\overline{U}(y)\leq\rho U(y)\text{.}
	\end{align}
	Then there exists a sequence $\epsilon^N\in\icoi$, $N\in\N$, with $\lim_{N\to\infty}\epsilon^N=0$ such that it holds for all $N\in\N$, $t\in\cc{0,T}$ that 
	\begin{equation}
		\begin{aligned}
			\E\!\lrcc{\exp\!\biggoo{\!e^{-\rho (\minofvttn)}U\oo{Y^N_t}\!+\!
					\int_{0}^{\minofv{t\!}{\!\tau^N}}
					\!\!\!\!e^{-\rho r}\overline{U}\oo{Y^N_r}\,dr\!}\!}
			\!\leq\E\bigcc{\exp\bigoo{U(Y^N_{0})}}e^{\epsilon^Nt}\text{.}
		\end{aligned}
	\end{equation}
\end{lemma}
\noindent
\textit{Proof of Lemma \ref{lemma 2}.} 
For all $N\in\N$ let $\lrfloor{\cdot}_{N}\colon\cc{0,T}\to\cc{0,T}$ be the function which satisfies for all $t\in\oc{0,T}$ that $\lrfloor{t}_{N}=\sup(\gk{0,\frac{T}{N},\frac{2T}{N},\dots,T}\cap\co{0,t})$ and $\lrfloor{0}_{N}=0$.
Choose $p\in\N$, $c\in\co{T^{\frac{1}{32}},\infty}$, and $N_0\in\N$ such that for all $x,y\in\Rd$, $N\in\N$ with $N\geq N_0$ it holds that 
\begin{align}
	\label{11}\norm{\mu(x)-\mu(y)}+\normf{\sigma(x)-\sigma(y)}\leq c\oo{1+\norm{x}^p+\norm{y}^p}\norm{x-y}\text{,}\\
	\label{12}\bb{\overline{U}(x)}+\normlrd{\Hess U(x)}+\norm{\nabla U(x)}+\bb{U(x)}+\norm{\mu(x)}+\normf{\sigma(x)}\leq c\oo{1+\norm{x}^p}\text{,}\\
	\label{13}e^{\sqrt{\bb{\log(N/T)}}}\leq c\bigoo{\tfrac{N}{T}}^{\frac{1}{32p}} \text{,}\\
	\label{14} c^p\bigoo{\tfrac{T}{N}}^{\frac{7}{32}}\bigoo{T^{\frac{3}{4}}+\sqrt{m}}\leq\bigoo{\tfrac{N}{T}}^{\frac{1}{32p}}\text{.}
\end{align}
For every $N\in\N$ let $\epsilon^N\colon\N\to\icoi$ satisfy that 
\begin{equation}
	\label{16}
	\begin{aligned}
		&\epsilon^N=\exp\biggoo{4^{p+\frac{1}{2}}c^{2p+3}\bigoo{\tfrac{T}{N}}^{\frac{3}{16}}\bigoo{T^{\frac{3}{4}}
				+\sqrt{m}}+\bigoo{c^2+3^pc^{2p+1}}\bigoo{\tfrac{T}{N}}^{\frac{31}{32}}}\\
		&\cdot\biggoo{4c^{2p+2}3^{2p}\bigoo{\tfrac{N}{T}}^{\frac{1}{16}}	\Bigoo{2c^{p+1}\bigoo{\tfrac{T}{N}}^{\frac{7}{32}}
				\bigoo{T^{\frac{3}{4}}+\sqrt{m}}+16\sqrt{m}\bigoo{\tfrac{T}{N}}^{\frace}}\\
			&+\Bigcc{104\sqrt{m}c^{p+1}\bigoo{\tfrac{T}{N}}^{\frac{15}{32}}
				\+ 2c^{2p+2}3^{p}\bigoo{\tfrac{T}{N}}^{\frac{3}{16}}\bigoo{T^{\frac{3}{4}}+\sqrt{m}}}\\
			&\hspace{2mm}\cdot\Bigcc{104\sqrt{m}c^{p+1}\bigoo{\tfrac{T}{N}}^{\frac{15}{32}}
				\+ 2c^{2p+2}3^{p}\bigoo{\tfrac{T}{N}}^{\frac{3}{16}}\bigoo{T^{\frac{3}{4}}\+\sqrt{m}}
				+4c^2\bigoo{\tfrac{N}{T}}^{\frac{1}{32}}}}
		\!\cdot\!\Bigoo{3^{2p}4c^{4p+2}\bigoo{\tfrac{N}{T}}^{\frac{1}{16}}}\text{.}
	\end{aligned}
\end{equation}
It is clear that $\lim_{N\to\infty}\epsilon^N=0$.
Lemma \ref{lemma 5} (applied for all $N\in\N$, $t\in\cc{0,T}$) and
Itô's formula impliey that for all $N\in\N$, $t\in\cc{0,T}$ it holds a.s. that
\begin{equation}
	\begin{aligned} 
		&\exp\biggoo{e^{-\rho \oo{\minofvttn}}U\oo{Y_t^N}+\int_0^{\minofvttn}e^{-\rho r}\overline{U}\oo{Y_r^N}\,dr}-e^{U(Y^N_{0})}\\
		=&\exp\biggoo{e^{-\rho \oo{\minofvttn}}U\oo{Y_{\minofvttn}^N}+\int_0^{\minofvttn}e^{-\rho r}\overline{U}\oo{Y_r^N}\,dr}-e^{U(Y^N_{0})}\\
		=&\int_0^{\minofvttn}\exp\biggoo{e^{-\rho s}U\oo{Y_s^N}+\int_{0}^s e^{-\rho r}\overline{U}\oo{Y_r^N}\,dr}
		e^{-\rho s}\bigoo{-\rho U\oo{Y_s^N}+\overline{U}\oo{Y_s^N}}\,ds\\
		&+\sum_{i=1}^d\int_0^{\minofvttn}\exp\biggoo{e^{-\rho s}U\oo{Y_s^N}+\int_{0}^s e^{-\rho r}\overline{U}\oo{Y_r^N}\,dr}e^{-\rho s}\Bigoo{\bigoo{\tfpartial{a_i}U}\oo{Y_s^N}}\,dY_s^N(i)\\
		&+\frace\sum_{i,k=1}^d\int_0^{\minofvttn}\exp\biggoo{e^{-\rho s}U\oo{Y_s^N}+\int_{0}^s e^{-\rho r}\overline{U}\oo{Y_r^N}\,dr}e^{-\rho s}\\
		&\hspace{5mm}\cdot\Bigoo{e^{-\rho s}\bigoo{\tfpartial{a_i}U}\oo{Y_s^N}\bigoo{\tfpartial{a_k}U}\oo{Y_s^N}+\bigoo{\tfrac{\partial^2}{\partial a_i\,\partial a_k}U}\oo{Y_s^N}}\,d[Y^N(i),Y^N(k)]_s\text{.}\\
	\end{aligned}
\end{equation}
This and (\ref{24}) show that for all $N\in\N$, $t\in\cc{0,T}$ it holds a.s. that
\begin{equation}
	\begin{aligned} 
		&\exp\biggoo{e^{-\rho \oo{\minofvttn}}U\oo{Y_t^N}+\int_0^{\minofvttn}e^{-\rho r}\overline{U}\oo{Y_r^N}\,dr}-e^{U(Y^N_{0})}\\	
		=&\int_0^{\minofvttn}\exp\biggoo{e^{-\rho s}U\oo{Y_s^N}+\int_{0}^s e^{-\rho r}\overline{U}\oo{Y_r^N}\,dr}e^{-\rho s}\bigoo{\overline{U}\oo{Y_s^N}-\rho U\oo{Y_s^N}}\,ds\\
		&+\sum_{i=1}^d\int_0^{\minofvttn}\exp\biggoo{e^{-\rho s}U\oo{Y_s^N}+\int_{0}^s e^{-\rho r}\overline{U}\oo{Y_r^N}\,dr}
		e^{-\rho s}\Bigoo{\bigoo{\tfpartial{a_i}U}\oo{Y_s^N}}\\
		&\hspace{5mm}\cdot \Bigoo{
				\mu(Y^N_{\lrfloor{s}_N})+\tfrace\sigma(Y^N_{\lrfloor{s}_N})\Delta\Pi^N(W_s-W_{\lrfloor{s}_N})
			\,ds+\sigma(Y^N_{\lrfloor{s}_N})D\Pi^N\oo{W_s\!-\!W_{\lrfloor{s}_N}}dW_s}_i\\
		&+\frace\sum_{i,k=1}^d\int_0^{\minofvttn}\exp\biggoo{e^{-\rho s}U\oo{Y_s^N}+\int_{0}^s e^{-\rho r}\overline{U}\oo{Y_r^N}\,dr}e^{-\rho s}\\
		&\hspace{5mm}\cdot\Bigoo{e^{-\rho s}\bigoo{\tfpartial{a_i}U}\oo{Y_s^N}\bigoo{\tfpartial{a_k}U}\oo{Y_s^N}+
			\bigoo{\tfrac{\partial^2}{\partial a_i\,\partial a_k}U}\oo{Y_s^N}}\\
		&\hspace{5mm}\cdot\bigoo{\sigma_i(Y^N_{\lrfloor{s}_N})D\Pi^N(W_s-W_{\lrfloor{s}_N})}\bigoo{\sigma_k(Y^N_{\lrfloor{s}_N})D\Pi^N(W_s-W_{\lrfloor{s}_N})}^{*}\,ds\text{.}
	\end{aligned}
\end{equation}
This implies that for all $N\in\N$, $t\in\cc{0,T}$ it holds a.s. that
\begin{equation}
	\label{43}
	\begin{aligned}
		&\exp\biggoo{e^{-\rho \oo{\minofvttn}}U\oo{Y_t^N}+\int_0^{\minofvttn}e^{-\rho r}\overline{U}\oo{Y_r^N}\,dr}-e^{U(Y^N_{0})}\\
		=\!&\!\int_0^{\minofvttn}\exp\biggoo{e^{-\rho s}U\oo{Y_s^N}+\int_{0}^s e^{-\rho r}\overline{U}\oo{Y_r^N}\,dr}e^{-\rho s}
		\cdot\biggoo{\overline{U}\oo{Y_s^N}-\rho U\oo{Y_s^N}\\
			&\hspace{5mm}+\biglr{\nabla U\oo{Y_s^N}}{\mu(Y^N_{\lrfloor{s}_N})+\tfrace\sigma(Y^N_{\lrfloor{s}_N})\Delta\Pi^N(W_s-W_{\lrfloor{s}_N})}\\
			&\hspace{5mm}+\tfrace
			\Biglr{\!\bigoo{\sigma(Y^N_{\lrfloor{s}_N})D\Pi^N(W_s-W_{\lrfloor{s}_N})}\!}
			{\!\Bigoo{e^{-\rho s}\bigoo{\nabla U\cdot (\nabla U)^*}(Y_s)+\!\Hess U\oo{Y_s^N}}\\
				&\hspace{14,5mm}\bigoo{\sigma(Y^N_{\lrfloor{s}_N})D\Pi^N(W_s-W_{\lrfloor{s}_N})}\!}_\tF
		}ds\\
		&+\int_0^{\minofvttn}\exp\biggoo{e^{-\rho s}U\oo{Y_s^N}+\int_{0}^s e^{-\rho r}\overline{U}\oo{Y_r^N}\,dr}e^{-\rho s}\\
		&\hspace{5mm}\cdot\biglr{\nabla U\oo{Y_s^N}}{\sigma(Y^N_{\lrfloor{s}_N})D\Pi^N\oo{W_s\!-\!W_{\lrfloor{s}_N}}dW_s}\text{.}
	\end{aligned}
\end{equation}
For all $n\in\N$ let $\delta_n\colon\Omega\to[0,T]$ be the function which satisfies that $\delta_n=\inf(\{s\in[0,T]\colon\norm{W_s}>n\}\cup\{T\})$.
The fact that $\overline{U}$, $U$, $\nabla U$, and $D\Pi^N$ are bounded on compact sets yields for every $n\in\N$ that the expectation of 
the stochastic integral on the righthandside of (\ref{43}) stopped at $\delta_n$ exists and vanishes.
This and (\ref{43}) show for all $N,n\in\N$, $t\in\cc{0,T}$ that
\begin{equation}
	\begin{aligned}
		&\E\lrcc{\exp\biggoo{e^{-\rho (\minminoft)}U\bigoo{Y^N_{\minofv{t}{\delta_n}}}+\int_0^{\minminoft} e^{-\rho r}\overline{U}\oo{Y_r^N}\,dr}}-e^{U(Y^N_{0})}\\
		=\,&\E\Bigg[\int_0^{\minminoft}\exp\biggoo{e^{-\rho s}U\oo{Y_s^N}+\int_{0}^s e^{-\rho r}\overline{U}\oo{Y_r^N}\,dr}e^{-\rho s}
		\cdot\biggoo{\overline{U}\oo{Y_s^N}-\rho U\oo{Y_s^N}\\
			&\hspace{4mm}+\biglr{\nabla U(Y^N_s)}{\mu(Y^N_s)}+\tfrace\biglr{\sigma(Y^N_s)}{\Hess U(Y^N_s)\sigma(Y^N_s)}_{\tF}\\
				&\hspace{10mm}+\tfrace e^{-\rho 	s}\bignorm{\oo{\sigma\oo{Y_s^N}}^*\nabla U\oo{Y_s^N}}^2\\
			&\hspace{4mm}+\lrr{\nabla U\oo{Y_s^N}}{\mu(Y^N_{\lrfloor{s}_N})-\mu\oo{Y_s^N}+\tfrace\sigma(Y^N_{\lrfloor{s}_N})\Delta\Pi^N(W_s-W_{\lrfloor{s}_N})}\\
			&\hspace{4mm}+\!\tfrace
			\Biglr{\!\bigoo{\sigma(Y^N_{\lrfloor{s}_N})D\Pi^N(W_s\- W_{\lrfloor{s}_N})}\!}
			{\!\!\bigoo{e^{-\rho s}\nabla U(\nabla U)^*\+\Hess U}\oo{Y_s^N}\\
				&\hspace{10mm}\cdot\bigoo{\sigma(Y^N_{\lrfloor{s}_N})D\Pi^N(W_s\- W_{\lrfloor{s}_N})}}_\tF\\
			&\hspace{4mm}-\!\tfrace
			\Biglr{\sigma(Y^N_s)}
			{\!\!\bigoo{e^{-\rho s}\nabla U(\nabla U)^*+\Hess U}\oo{Y_s^N}\sigma(Y^N_s)}_\tF
		}ds\Bigg]\!\text{.}
	\end{aligned}
\end{equation}
This implies for all $N,n\in\N$, $t\in\cc{0,T}$ that
\begin{equation}
	\begin{aligned}
		&\E\lrcc{\exp\biggoo{e^{-\rho (\minminoft)}U\bigoo{Y^N_{\minofv{t}{\delta_n}}}+\int_0^{\minminoft} e^{-\rho r}\overline{U}\oo{Y_r^N}\,dr}}\\
		&-\E\lrcc{\exp\biggoo{e^{-\rho (\minminoftab)}U\bigoo{Y^N_{\minofv{\lrfloor{t}_N}{\delta_n}}}+\int_0^{\minminoftab} e^{-\rho r}\overline{U}\oo{Y_r^N}\,dr}}\\
		=\,&\E\Bigg[\int_{\minminoftab}^{\minminoft}\exp\biggoo{e^{-\rho s}U\oo{Y_s^N}+\int_{0}^s e^{-\rho r}\overline{U}\oo{Y_r^N}\,dr}e^{-\rho s}
		\cdot\biggoo{\overline{U}\oo{Y_s^N}-\rho U\oo{Y_s^N}\\
			&\hspace{4mm}+\biglr{\nabla U(Y^N_s)}{\mu(Y^N_s)}+\tfrace\biglr{\sigma(Y^N_s)}{\Hess U(Y^N_s)\sigma(Y^N_s)}_{\tF}\\
				&\hspace{10mm}+\tfrace e^{-\rho 	s}\bignorm{\oo{\sigma\oo{Y_s^N}}^*\nabla U\oo{Y_s^N}}^2\\
			&\hspace{4mm}+\lrr{\nabla U\oo{Y_s^N}}{\mu(Y^N_{\lrfloor{s}_N})-\mu\oo{Y_s^N}+\tfrace\sigma(Y^N_{\lrfloor{s}_N})\Delta\Pi^N(W_s-W_{\lrfloor{s}_N})}\\
			&\hspace{4mm}+\!\tfrace
			\Biglr{\!\bigoo{\sigma(Y^N_{\lrfloor{s}_N})D\Pi^N(W_s\- W_{\lrfloor{s}_N})}\!}
			{\!\!\bigoo{e^{-\rho s}\nabla U(\nabla U)^*\+ \Hess U}\oo{Y_s^N}\\
				&\hspace{10mm}\cdot\bigoo{\sigma(Y^N_{\lrfloor{s}_N})D\Pi^N(W_s\- W_{\lrfloor{s}_N})}}_\tF\\
			&\hspace{4mm}-\!\tfrace
			\Biglr{\sigma(Y^N_s)}
			{\!\!\bigoo{e^{-\rho s}\nabla U(\nabla U)^*+\Hess U}\oo{Y_s^N}\sigma(Y^N_s)}_\tF
		}ds\Bigg]\!\text{.}
	\end{aligned}
\end{equation}
Assumption (\ref{1}), $p\in\icoi$, the Cauchy-Schwarz inequality, and the fact that for all symmetric $A\in\Rdd$ and all $y,z\in\Rdm$ it holds that $\lrrf{y}{Ay}-\lrrf{z}{Az}=\lrrf{y+z}{A(y-z)}\leq\norm{A}_{L(\Rd)}\,\normf{z-y}\,\normf{z+y}$ imply that for all $N,n\in\N$, $t\in\cc{0,T}$ it holds that
\begin{equation}
	\begin{aligned}
		&\E\lrcc{\exp\biggoo{e^{-\rho (\minminoft)}U\bigoo{Y^N_{\minofv{t}{\delta_n}}}+\int_0^{\minminoft} e^{-\rho r}\overline{U}\oo{Y_r^N}\,dr}}\\
		&-\E\lrcc{\exp\biggoo{e^{-\rho (\minminoftab)}U\bigoo{Y^N_{\minofv{\lrfloor{t}_N}{\delta_n}}}+\int_0^{\minminoftab} e^{-\rho r}\overline{U}\oo{Y_r^N}\,dr}}\\
		\leq\,&\E\Bigg[\int_{\minminoftab}^{\minminoft}\exp\biggoo{e^{-\rho s}U\oo{Y_s^N}+\int_{0}^s e^{-\rho r}\overline{U}\oo{Y_r^N}\,dr}e^{-\rho s}\\
		&\hspace{3mm}\cdot\biggoo{\Bignorm{\nabla U\oo{Y_s^N}}\;
			\Bignorm{\mu(Y^N_{\lrfloor{s}_N})-\mu\oo{Y_s^N}+\tfrace\sigma(Y^N_{\lrfloor{s}_N})\Delta\Pi^N(W_s-W_{\lrfloor{s}_N})}\\
			&\hspace{5mm}\+\tfrace	\Bignormlrd{\bigoo{e^{-\rho s}\nabla U(\nabla U)^*\!+\!\Hess U}\!\oo{Y_s^N}}
			\Bignormf{\sigma(Y^N_s)\!-\!\sigma(Y^N_{\lrfloor{s}_N})D\Pi^N\!(W_s\!-\!W_{\lrfloor{s}_N})}\\
			&\hspace{9,5mm}\cdot\Bignormf{\sigma(Y^N_s)+\sigma(Y^N_{\lrfloor{s}_N})D\Pi^N(W_s-W_{\lrfloor{s}_N})}
		}ds\Bigg]\!\text{.}
	\end{aligned}
\end{equation}
The triangle inequality and the fact that $\forall x\in\Rd\colon\normlrd{xx^*}\leq\norm{x}^2$ show that for all $N,n\in\N$, $t\in\cc{0,T}$ it holds that 
\begin{equation}
	\begin{aligned}
		\label{25}
		&\E\lrcc{\exp\biggoo{e^{-\rho (\minminoft)}U\bigoo{Y^N_{\minofv{t}{\delta_n}}}+\int_0^{\minminoft} e^{-\rho r}\overline{U}\oo{Y_r^N}\,dr}}\\
		&-\E\lrcc{\exp\biggoo{e^{-\rho (\minminoftab)}U\bigoo{Y^N_{\minofv{\lrfloor{t}_N}{\delta_n}}}+\int_0^{\minminoftab} e^{-\rho r}\overline{U}\oo{Y_r^N}\,dr}}\\
		\leq\,&\E\Bigg[\int_{\minminoftab}^{\minminoft}\exp\biggoo{e^{-\rho s}U\oo{Y_s^N}+\int_{0}^s e^{-\rho r}\overline{U}\oo{Y_r^N}\,dr}e^{-\rho s}\\
		&\hspace{10mm}\cdot\biggoo{\bignorm{\nabla U\oo{Y_s^N}}\;
			\Bigoo{\bignorm{\mu\oo{Y_s^N}-\mu(Y^N_{\lrfloor{s}_N})}+\tfrace\bignorm{\sigma(Y^N_{\lrfloor{s}_N})\Delta\Pi^N(W_s-W_{\lrfloor{s}_N})}}\\
			&\hspace{15mm}+\tfrace\Bigoo{\bignorm{\nabla U(Y_s^N)}^2+\bignormlrd{\Hess U\oo{Y_s^N}}}\\
			&\hspace{19,5mm}\cdot\Bigoo{\,\bignormf{\sigma(Y^N_{\lrfloor{s}_N})\cdot\bigoo{I_{\Rm}-D\Pi^N(W_s-W_{\lrfloor{s}_N})}}
				+\bignormf{\sigma(Y^N_s)-\sigma(Y^N_{\lrfloor{s}_N})}}\\
			&\hspace{19,5mm}\cdot\Bigoo{\,\bignormf{\sigma(Y^N_{\lrfloor{s}_N})\cdot\bigoo{I_{\Rm}-D\Pi^N(W_s-W_{\lrfloor{s}_N})}}
				+\bignormf{\sigma(Y^N_s)-\sigma(Y^N_{\lrfloor{s}_N})}\\
				&\hspace{24,5mm}+2\bignormf{\sigma(Y^N_{\lrfloor{s}_N})}}
		}ds\Bigg]\!\text{.}
	\end{aligned}
\end{equation}
The fundamental theorem of calculus, \ref{021} in Lemma \ref{lemma 6} (applied for every $N\in\N\cap\co{N_0,\infty}$, $t\in\cc{0,T}$), (\ref{12}), (\ref{13}), and (\ref{14}) show that for all $t\in\cc{0,T}$, $N\in\N$ with $N\geq N_0$ it holds a.s. that 
\begin{equation}
	\begin{aligned}
		\label{17}
		&\mathbbm{1}_{\norm{Y^N_{\lrfloor{t}_N}}\leq e^{\sqrt{\bb{\log(N/T)}}}}
		\cdot \Bigoo{U\bigoo{Y^N_t}-U\bigoo{Y^N_{\lrfloor{t}_N}}}\\
		=\,&\mathbbm{1}_{\norm{Y^N_{\lrfloor{t}_N}}\leq e^{\sqrt{\bb{\log(N/T)}}}}
		\cdot\int_0^1(\nabla U)^*\bigoo{Y^N_t\lambda+Y^N_{\lrfloor{t}_N}(1-\lambda)}d\lambda\bigoo{Y_t^N-Y^N_{\lrfloor{t}_N}}\\
		\leq\,& c\mathbbm{1}_{\norm{Y^N_{\lrfloor{t}_N}}\leq e^{\sqrt{\bb{\log(N/T)}}}}
		\cdot\int_0^1\bigoo{1+\bignorm{Y_t^N\lambda+Y^N_{\lrfloor{t}_N}(1-\lambda)}^p}d\lambda\bignorm{Y_t^N-Y^N_{\lrfloor{t}_N}}\\
		\leq\,& c\mathbbm{1}_{\norm{Y^N_{\lrfloor{t}_N}}\leq e^{\sqrt{\bb{\log(N/T)}}}}
		\cdot\Bigoo{1+\bignorm{Y_t^N-Y^N_{\lrfloor{t}_N}}+\bignorm{Y^N_{\lrfloor{t}_N}}}^p\bignorm{Y_t^N-Y^N_{\lrfloor{t}_N}}\\
		\leq\,& c\Bigoo{2c^{p+1}\bigoo{\tfrac{T}{N}}^{\frac{7}{32}}\bigoo{T^{\frac{3}{4}}+\sqrt{m}}+2e^{\sqrt{\bb{\log(N/T)}}}}^p
		2c^2\bigoo{\tfrac{T}{N}}^{\frac{7}{32}}\bigoo{T^{\frac{3}{4}}+\sqrt{m}}\\
		\leq\,& \Bigoo{2c^{p+1}\bigoo{\tfrac{T}{N}}^{\frac{7}{32}}\bigoo{T^{\frac{3}{4}}+\sqrt{m}}+2c\bigoo{\tfrac{N}{T}}^{\frac{1}{32p}}}^p
		2c^3\bigoo{\tfrac{T}{N}}^{\frac{7}{32}}\bigoo{T^{\frac{3}{4}}+\sqrt{m}}\\
		\leq\,& 4^pc^{p}\bigoo{\tfrac{N}{T}}^{\frac{1}{32}}
		2c^3\bigoo{\tfrac{T}{N}}^{\frac{7}{32}}\bigoo{T^{\frac{3}{4}}+\sqrt{m}}\\
		=\,& 4^{p+\frac{1}{2}}c^{p+3}\bigoo{\tfrac{T}{N}}^{\frac{3}{16}}\bigoo{T^{\frac{3}{4}}+\sqrt{m}}\text{.}
	\end{aligned}
\end{equation}
Next \ref{022} in lemma \ref{lemma 6} (applied for every $N\in\N\cap\co{N_0,\infty}$, $t\in\cc{0,T}$) implies that for all $t\in\cc{0,T}$, $N\in\N$ with $N\geq N_0$ it holds a.s. that
\begin{equation}
	\begin{aligned}
		\label{29}
		&\mathbbm{1}_{\norm{Y^N_{\lrfloor{t}_N}}\leq e^{\sqrt{\bb{\log(N/T)}}}}\cdot
		\tfrace\Bigoo{c^2\bigoo{1+\norm{Y_t^N}^p}^2+c\bigoo{1+\norm{Y_t^N}^p}}\\
		\leq\,&\mathbbm{1}_{\norm{Y^N_{\lrfloor{t}_N}}\leq e^{\sqrt{\bb{\log(N/T)}}}}\cdot c^2\bigoo{1+\norm{Y_t^N}^p}^2\\
		\leq\,&c^2\Bigoo{1+3^pc^{p}\bigoo{\tfrac{N}{T}}^{\frac{1}{32}}}^2\\
		\leq\,&3^{2p}4c^{2p+2}\bigoo{\tfrac{N}{T}}^{\frac{1}{16}}\text{.}
	\end{aligned}
\end{equation}
Moreover (\ref{12}) and (\ref{13}) show for all $t\in\cc{0,T}$, $N\in\N$ with $N\geq N_0$ that
\begin{equation}
	\begin{aligned}
		\label{30}
		&\mathbbm{1}_{\norm{Y^N_{\lrfloor{t}_N}}\leq e^{\sqrt{\bb{\log(N/T)}}}}
		\cdot\bigoo{e^{-\rho t}-e^{-\rho\lrfloor{t}_N}}U(Y_{\lrfloor{s}_N}^N)\\
		\leq\,&\mathbbm{1}_{\norm{Y^N_{\lrfloor{t}_N}}\leq e^{\sqrt{\bb{\log(N/T)}}}}
		\cdot\rho(t-\lrfloor{t}_N)c(1+\norm{Y^N_{\lrfloor{t}_N}}^p)\\
		\leq\,&\rho\bigoo{\tfrac{T}{N}}2ce^{p\sqrt{\bb{\log(N/T)}}}\\
		\leq\,&2\rho c^{p+1}\bigoo{\tfrac{T}{N}}^{\frac{31}{32}}\text{.}
	\end{aligned}
\end{equation}
Next (\ref{12}), \ref{022} in Lemma \ref{lemma 6} (applied for every $N\in\N\cap\co{N_0,\infty}$, $t\in\cc{0,T}$), and the fact that $\sup_{x\in\cc{0,\minofvttn}}e^{-\rho x}\leq 1$ imply that for all $t\in\cc{0,T}$, $N\in\N$ with $N\geq N_0$ it holds a.s. that
\begin{equation}
	\begin{aligned}
		\label{37}
		&\int_{\minofvtabtn}^{\minofvttn}e^{-\rho r}\overline{U}(Y_r^N)dr\\
		\leq\,&\int_{\minofvtabtn}^{\minofvttn}e^{-\rho r}c\bigoo{1+\norm{Y_r^N}^p}dr\\
		\leq\,&c\bigoo{1+3^pc^{p}\bigoo{\tfrac{N}{T}}^{\frac{1}{32}}}\bigoo{\minofv{t\!\!}{\!\!\tau_N}-\minofv{\lrfloor{t}_N\!\!}{\!\!\tau_N}}\\
		\leq\,&\bigoo{c^2+3^pc^{p+1}}\bigoo{\tfrac{T}{N}}^{\frac{31}{32}}\text{.}
	\end{aligned}
\end{equation}
Inserting \ref{024} and \ref{025} from lemma \ref{lemma 6} (applied for every $N\in\N\cap\co{N_0,\infty}$, $t\in\cc{0,T}$), and (\ref{12}) in (\ref{25}) show that for all $t\in\cc{0,T}$, $N,n\in\N$ with $N\geq N_0$ it holds a.s. that 
\begin{equation}
	\begin{aligned}
		&\E\lrcc{\exp\biggoo{e^{-\rho (\minminoft)}U\bigoo{Y^N_{\minofv{t}{\delta_n}}}+\int_0^{\minminoft} e^{-\rho r}\overline{U}\oo{Y_r^N}\,dr}}\\
		&-\E\lrcc{\exp\biggoo{e^{-\rho (\minminoftab)}U\bigoo{Y^N_{\minofv{\lrfloor{t}_N}{\delta_n}}}+\int_0^{\minminoftab} e^{-\rho r}\overline{U}\oo{Y_r^N}\,dr}}\\
		\leq\,&\E\Biggcc{\exp\biggoo{e^{-\rho(\minminoftab)}\,U(Y^N_{\minofv{\lrfloor{t}_N}{\delta_n}})+\int_{0}^{\minminoftab}e^{-\rho r}\overline{U}(Y_r^N)dr}\\
			&\cdot\!\!\int_{\minminoftab}^{\minminoft}\!\!\!\!\!\!\!\!\exp\!\biggoo{\!e^{-\rho s}\bigoo{U(Y_s^N)\-U(Y_{\lrfloor{s}_N}^N)}
				\+\bigoo{e^{-\rho s}\-e^{-\rho\lrfloor{s}_N}}U(Y_{\lrfloor{s}_N}^N)
				\+\!\int_{\minofv{\lrfloor{s}_N\!}{\!\tau^N}}^{\minofv{s\!}{\!\tau^N}}\!\!\!\!\!\!\!\!e^{-\rho r}\overline{U}(Y_r^N)dr}\\
			&\hspace{10mm}\cdot\biggoo{4c^{2p+2}3^{2p}\bigoo{\tfrac{N}{T}}^{\frac{1}{16}}	\Bigoo{2c^{p+1}\bigoo{\tfrac{T}{N}}^{\frac{7}{32}}
					\bigoo{T^{\frac{3}{4}}+\sqrt{m}}+\tfrace\bignorm{\Delta\Pi^N(W_s-W_{\lrfloor{s}_N})}_{L(\Rm)}}\\
				&\hspace{15mm}+\tfrace\Bigoo{c^2\bigoo{1+\norm{Y_s^N}^p}^2+c\bigoo{1+\norm{Y_s^N}^p}}\\
				&\hspace{19,5mm}\cdot\Bigoo{\bignorm{D\Pi^N(W_s\!-\!W_{\lrfloor{s}_N})\!-\!I_{\Rm}}_{L(\Rm)}
					\cdot 2c^{p+1}\bigoo{\tfrac{N}{T}}^{\frac{1}{32}}
					+2c^{2p+2}3^{p}\bigoo{\tfrac{T}{N}}^{\frac{3}{16}}\bigoo{T^{\frac{3}{4}}\!+\!\sqrt{m}}}\\
				&\hspace{19,5mm}\cdot\Bigoo{\bignorm{D\Pi^N(W_s\!-\!W_{\lrfloor{s}_N})\!-\!I_{\Rm}}_{L(\Rm)}
					\cdot 2c^{p+1}\bigoo{\tfrac{N}{T}}^{\frac{1}{32}}
					+2c^{2p+2}3^{p}\bigoo{\tfrac{T}{N}}^{\frac{3}{16}}\bigoo{T^{\frac{3}{4}}\!+\!\sqrt{m}}\\
					&\hspace{24,5mm}+\!2c\bigoo{1\!+\!\norm{Y_{\lrfloor{s}_N}^N}^p}}
			}ds}\!\text{.}
	\end{aligned}
\end{equation}
Moreover (\ref{17}), (\ref{29}), (\ref{30}), and (\ref{37}) imply that for all $t\in\cc{0,T}$, $N,n\in\N$ with $N\geq N_0$ it holds a.s. that 
\begin{equation}
	\begin{aligned}
		&\E\lrcc{\exp\biggoo{e^{-\rho (\minminoft)}U\bigoo{Y^N_{\minofv{t}{\delta_n}}}+\int_0^{\minminoft} e^{-\rho r}\overline{U}\oo{Y_r^N}\,dr}}\\
		&-\E\lrcc{\exp\biggoo{e^{-\rho (\minminoftab)}U\bigoo{Y^N_{\minofv{\lrfloor{t}_N}{\delta_n}}}+\int_0^{\minminoftab} e^{-\rho r}\overline{U}\oo{Y_r^N}\,dr}}\\
		\leq\,&\E\Biggcc{\exp\biggoo{e^{-\rho(\minminoftab)}\,U(Y^N_{\minofv{\lrfloor{t}_N}{\delta_n}})+\int_{0}^{\minminoftab}e^{-\rho r}\overline{U}(Y_r^N)dr}\\
			&\cdot\int_{\minminoftab}^{\minminoft}\!\!\!\!\!\!\!\!\exp\biggoo{e^{-\rho s}4^{p+\frac{1}{2}}c^{2p+3}\bigoo{\tfrac{T}{N}}^{\frac{3}{16}}\bigoo{T^{\frac{3}{4}}\+\sqrt{m}}
				\+2\rho c^{p+1}\bigoo{\tfrac{T}{N}}^{\frac{31}{32}}
				\+\bigoo{c^2+3^pc^{p+1}}\bigoo{\tfrac{T}{N}}^{\frac{31}{32}}}\\
			&\hspace{10mm}\cdot\biggoo{4c^{2p+2}3^{2p}\bigoo{\tfrac{N}{T}}^{\frac{1}{16}}\Bigoo{2c^{p+1}\bigoo{\tfrac{T}{N}}^{\frac{7}{32}}\bigoo{T^{\frac{3}{4}}+\sqrt{m}}
					+\tfrace\bignorm{\Delta\Pi^N(W_s-W_{\lrfloor{s}_N})}_{L(\Rm)}}\\
				&\hspace{15mm}+\Bigoo{\bignorm{D\Pi^N(W_s\!-\!W_{\lrfloor{s}_N})\!-\!I_{\Rm}}_{L(\Rm)}\cdot 2c^{p+1}\bigoo{\tfrac{N}{T}}^{\frac{1}{32}}
					+2c^{2p+2}3^{p}\bigoo{\tfrac{T}{N}}^{\frac{3}{16}}\bigoo{T^{\frac{3}{4}}\!+\!\sqrt{m}}}\\
				&\hspace{19,5mm}\cdot\Bigoo{\bignorm{D\Pi^N(W_s\!-\!W_{\lrfloor{s}_N})\!-\!I_{\Rm}}_{L(\Rm)}\cdot 2c^{p+1}\bigoo{\tfrac{N}{T}}^{\frac{1}{32}}
					+2c^{2p+2}3^{p}\bigoo{\tfrac{T}{N}}^{\frac{3}{16}}\bigoo{T^{\frac{3}{4}}\!+\!\sqrt{m}}\\
					&\hspace{24,5mm}+\!4ce^{p\sqrt{\bb{\log(N/T)}}}}}
			\cdot\Bigoo{3^{2p}4c^{2p+2}\bigoo{\tfrac{N}{T}}^{\frac{1}{16}}
			}ds}\!\text{.}
	\end{aligned}
\end{equation}
Hence the fact that for all $t\in\cc{0,T},N\in\N$ it holds that $Y_{\lrfloor{t}_N}^N$ and $(W_s-W_{\lrfloor{t}_N})_{s\in\cc{\lrfloor{t}_N,t}}$ are independent, Fubini's theorem, Hölder's inequality, and (\ref{13}) imply that for all $t\in\cc{0,T}$, $N,n\in\N$ with $N\geq N_0$ it holds a.s. that
\begin{equation}
	\begin{aligned}
		&\E\lrcc{\exp\biggoo{e^{-\rho (\minminoft)}U\bigoo{Y^N_{\minofv{t}{\delta_n}}}+\int_0^{\minminoft} e^{-\rho r}\overline{U}\oo{Y_r^N}\,dr}}\\
		&-\E\lrcc{\exp\biggoo{e^{-\rho (\minminoftab)}U\bigoo{Y^N_{\minofv{\lrfloor{t}_N}{\delta_n}}}+\int_0^{\minminoftab} e^{-\rho r}\overline{U}\oo{Y_r^N}\,dr}}\\
		\leq\,&\E\Biggcc{\exp\biggoo{e^{-\rho(\minminoftab)}\,U(Y^N_{\minofv{\lrfloor{t}_N}{\delta_n}})+\int_{0}^{\minminoftab}e^{-\rho r}\overline{U}(Y_r^N)dr}}\\
		&\cdot\int_{\lrfloor{t}_N}^{t}\exp\biggoo{e^{-\rho s}4^{p+\frac{1}{2}}c^{2p+3}\bigoo{\tfrac{T}{N}}^{\frac{3}{16}}\bigoo{T^{\frac{3}{4}}\+\sqrt{m}}
			\+2\rho c^{p+1}\bigoo{\tfrac{T}{N}}^{\frac{31}{32}}
			\+\bigoo{c^2+3^pc^{p+1}}\bigoo{\tfrac{T}{N}}^{\frac{31}{32}}}\\
		&\hspace{5mm}\cdot\biggoo{4c^{2p+2}3^{2p}\bigoo{\tfrac{N}{T}}^{\frac{1}{16}}	\Bigoo{2c^{p+1}\bigoo{\tfrac{T}{N}}^{\frac{7}{32}}
				\bigoo{T^{\frac{3}{4}}+\sqrt{m}}+\tfrace\bignorm{\Delta\Pi^N(W_s-W_{\lrfloor{s}_N})}_{L^2(\Omega;L(\Rm))}}\\
			&\hspace{10mm}+\Bigoo{\bignorm{D\Pi^N(W_s\!-\!W_{\lrfloor{s}_N})\!-\!I_{\Rm}}_{\LzweiOLRm}\!\cdot\! 2c^{p+1}\bigoo{\tfrac{N}{T}}^{\frac{1}{32}}
				\+ 2c^{2p+2}3^{p}\bigoo{\tfrac{T}{N}}^{\frac{3}{16}}\bigoo{T^{\frac{3}{4}}\!+\!\sqrt{m}}\!}\\
			&\hspace{12mm}\cdot\Bigoo{\bignorm{D\Pi^N(W_s\!-\!W_{\lrfloor{s}_N})\!-\!I_{\Rm}}_{\LzweiOLRm}\!\cdot\! 2c^{p+1}\bigoo{\tfrac{N}{T}}^{\frac{1}{32}}
				\+ 2c^{2p+2}3^{p}\bigoo{\tfrac{T}{N}}^{\frac{3}{16}}\bigoo{T^{\frac{3}{4}}\!+\!\sqrt{m}}\\
				&\hspace{19,5mm}+\!4c^2\bigoo{\tfrac{N}{T}}^{\frac{1}{32}}}\!}
		\cdot\Bigoo{3^{2p}4c^{4p+2}\bigoo{\tfrac{N}{T}}^{\frac{1}{16}}
		}ds\text{.}
	\end{aligned}
\end{equation}
Lemma \ref{lemma 1} (applied for every $N\in\N\cap\co{N_0,\infty}$, $s\in\cc{0,T}$ with $h\curvearrowleft T/N$, $W\curvearrowleft (W_{\lrfloor{s}_N+t}-W_{\lrfloor{s}_N})_{t\in\cc{0,T/N}}$) shows that for all $t\in\cc{0,T}$, $N,n\in\N$ with $N\geq N_0$ it holds a.s. that
\begin{equation}
	\begin{aligned}
		&\E\lrcc{\exp\biggoo{e^{-\rho (\minminoft)}U\bigoo{Y^N_{\minofv{t}{\delta_n}}}+\int_0^{\minminoft} e^{-\rho r}\overline{U}\oo{Y_r^N}\,dr}}\\
		&-\E\lrcc{\exp\biggoo{e^{-\rho (\minminoftab)}U\bigoo{Y^N_{\minofv{\lrfloor{t}_N}{\delta_n}}}+\int_0^{\minminoftab} e^{-\rho r}\overline{U}\oo{Y_r^N}\,dr}}\\
		\leq\,&\E\Biggcc{\exp\biggoo{e^{-\rho(\minminoftab)}\,U(Y^N_{\minofv{\lrfloor{t}_N}{\delta_n}})+\int_{0}^{\minminoftab}e^{-\rho r}\overline{U}(Y_r^N)dr}}\\
		&\cdot\exp\biggoo{4^{p+\frac{1}{2}}c^{2p+3}\bigoo{\tfrac{T}{N}}^{\frac{3}{16}}\bigoo{T^{\frac{3}{4}}
				+\sqrt{m}}+\bigoo{c^2+3^pc^{2p+1}}\bigoo{\tfrac{T}{N}}^{\frac{31}{32}}}\\
		&\cdot\biggoo{4c^{2p+2}3^{2p}\bigoo{\tfrac{N}{T}}^{\frac{1}{16}}	\Bigoo{2c^{p+1}\bigoo{\tfrac{T}{N}}^{\frac{7}{32}}
				\bigoo{T^{\frac{3}{4}}+\sqrt{m}}+16\sqrt{m}\bigoo{\tfrac{T}{N}}^{\frace}}\\
			&+\Bigcc{104\sqrt{m}c^{p+1}\bigoo{\tfrac{T}{N}}^{\frac{15}{32}}
				\+ 2c^{2p+2}3^{p}\bigoo{\tfrac{T}{N}}^{\frac{3}{16}}\bigoo{T^{\frac{3}{4}}+\sqrt{m}}}\\
			&\hspace{2mm}\cdot\Bigcc{104\sqrt{m}c^{p+1}\bigoo{\tfrac{T}{N}}^{\frac{15}{32}}
				\+ 2c^{2p+2}3^{p}\bigoo{\tfrac{T}{N}}^{\frac{3}{16}}\bigoo{T^{\frac{3}{4}}\+\sqrt{m}}
				+4c^2\bigoo{\tfrac{N}{T}}^{\frac{1}{32}}}}\\
		&\hspace{5mm}\!\cdot\!\Bigoo{3^{2p}4c^{4p+2}\bigoo{\tfrac{N}{T}}^{\frac{1}{16}}}\cdot\bigoo{t-\lrfloor{t}_N}\text{.}
	\end{aligned}
\end{equation}
This and (\ref{16}) show that for all $t\in\cc{0,T}$, $N,n\in\N$ with $N\geq N_0$ it holds a.s. that 
\begin{equation}
	\begin{aligned}
		&\E\lrcc{\exp\biggoo{e^{-\rho (\minminoft)}U\bigoo{Y^N_{\minofv{t}{\delta_n}}}+\int_0^{\minminoft} e^{-\rho r}\overline{U}\oo{Y_r^N}\,dr}}\\
		\leq\,&\E\Biggcc{\exp\biggoo{e^{-\rho(\minminoftab)}\,U(Y^N_{\minofv{\lrfloor{t}_N}{\delta_n}})\+\int_{0}^{\minminoftab}\!\!\!\!\!\!\!\!e^{-\rho r}\overline{U}(Y_r^N)dr}}
		\!\cdot\!\Bigoo{1\+\epsilon^N\bigoo{t\-\lrfloor{t}_N}}\\
		\leq\,&\E\Biggcc{\exp\biggoo{e^{-\rho(\minminoftab)}\,U(Y^N_{\minofv{\lrfloor{t}_N}{\delta_n}})\+\int_{0}^{\minminoftab}\!\!\!\!\!\!\!\!e^{-\rho r}\overline{U}(Y_r^N)dr}}
		\!\cdot\!\exp\!\Bigoo{\epsilon^N\!\bigoo{t\-\lrfloor{t}_N}}\text{.}
	\end{aligned}
\end{equation}
Iterating this inequality implies that for all $t\in\cc{0,T}$, $N,n\in\N$ with $N\geq N_0$ it holds a.s. that 
\begin{equation}
	\begin{aligned}
		&\E\lrcc{\exp\biggoo{e^{-\rho (\minminoft)}U\bigoo{Y^N_{\minofv{t}{\delta_n}}}+\int_0^{\minminoft} e^{-\rho r}\overline{U}\oo{Y_r^N}\,dr}}\\
		\leq\,&\E\bigcc{\exp\bigoo{U(Y^N_{0})}}\cdot e^{\epsilon^Nt}\text{.}
	\end{aligned}
\end{equation}
This, Fatous's Lemma, $\mathbb{P}(\lim_{n\to\infty}\delta_n=T)=1$, and path continuity show that for all $t\in\cc{0,T}$, $N,n\in\N$ with $N\geq N_0$ it holds a.s. that 
\begin{equation}
	\begin{aligned}
		&\E\!\lrcc{\exp\!\biggoo{\!e^{-\rho (\minofvttn)}U\oo{Y^N_t}\!+\!
				\int_{0}^{\minofv{t\!}{\!\tau^N}}
				\!\!\!\!e^{-\rho r}\overline{U}\oo{Y^N_r}\,dr\!}\!}\\
		\leq\,&\liminf_{n\to\infty}\E\lrcc{\exp\biggoo{e^{-\rho (\minminoft)}U\bigoo{Y^N_{\minofv{t}{\delta_n}}}+\int_0^{\minminoft} e^{-\rho r}\overline{U}\oo{Y_r^N}\,dr}}\\
		\leq\,&\E\bigcc{\exp\bigoo{U(Y^N_{0})}}\cdot e^{\epsilon^Nt}\text{.}
	\end{aligned}
\end{equation}
The proof of Lemma \ref{lemma 2} is thus completed. $ \hfill\square $\\
\section{Moment estimates}
The following lemma proves that stopped Brownian-increment tamed Euler approximations
have uniformly bounded moments.
Lemma \ref{lemma 3} is motivated by \cite{HutzenthalerJentzen2015} which proves the
respective result for increment-tamed Euler approximations.
\begin{lemma}
	\label{lemma 3}
	
	Let $\rho\in \icoi$,
	$T\in\iooi$,
	$d, m, p\in\N$, 
	$N\in\N\cap\co{T,\infty}$,
	$c\in\co{T^{\frac{1}{32}},\infty}$,
	let $\mu\colon\Rd\to\Rd$, 
	$\sigma\colon\Rd\to\Rdm$,
	$U\in C^2(\Rd,\co{0,\infty})$, 
	let $\Pi^N\in C^2(\Rm,\Rm)$ be the function which satisfies for all $x=\oo{x_1,\dots,x_m}\in\Rm$ that
	\begin{equation}
		\Pi^N(x)=\lroo{x_i\,\exp\lroo{-\frac{Nx_i^4}{T}}}_{i\in\gk{1,\dots,m}}\text{,}
	\end{equation}
	let $(\Omega,\F,\P,(\mathbb{F}_t)_{t\in\cc{0,T}})$ be a filtered probability space,
	let $W\colon\icc{0}{T}\times\Omega\to\R^m$ be a standard $(\mathbb{F}_t)_{t\in\cc{0,T}}$-Brownian motion with continuous sample paths,
	let $Y^N\colon\cc{0,T}\times\Omega\to\Rd$ satisfy that $Y^N_0$ is $\mathbb{F}_0$-measurable
	and for all $k\in\gk{0,1,\dots,N-1}$, $t\in\bigcc{\frac{kT}{N},\frac{(k+1)T}{N}}$ that
	\begin{align}
		\label{42}
		Y_t^N=Y^N_{\frac{kT}{N}}+\mathbbm{1}_{\norm{Y^N_{\frac{kT}{N}}}\leq e^{\sqrt{\bb{\log(N/T)}}}}
		\cdot\Bigoo{\mu\bigoo{Y^N_{\frac{kT}{N}}}\bigoo{t-\tfrac{kT}{N}}
			+\sigma\bigoo{Y^N_{\frac{kT}{N}}}\Pi^N\bigoo{W_t-W_{\frac{kT}{N}}}}\text{,}
	\end{align}
	and assume for all $x,y,z\in\Rd$ that  
	\begin{align}
		\label{34}	\lrr{z}{\bigoo{\Hess U(y)-\Hess U(x)}z}\leq\,c\norm{z}^2\bigoo{1+U(x)+\norm{y-x}^p}\norm{y-x}\text{,}\\
		\label{112}\normlrd{\Hess U(x)}+\norm{\nabla U(x)}+\bb{U(x)}+\norm{\mu(x)}+\normf{\sigma(x)}\leq c\oo{1+\norm{x}^p}\text{,}\\
		\label{113}e^{\sqrt{\bb{\log(N/T)}}}\leq c\bigoo{\tfrac{N}{T}}^{\frac{1}{32p}} \text{,}\\
		\label{114} c^p\bigoo{\tfrac{T}{N}}^{\frac{7}{32}}\bigoo{T^{\frac{3}{4}}+\sqrt{m}}\leq\bigoo{\tfrac{N}{T}}^{\frac{1}{32p}}\text{,}\\
		\label{33}	\!\!\!\!\!\!\lrr{\nabla U(y)}{\mu(y)}+\tfrace\lrrf{\sigma(y)}{\Hess U(y)\sigma(y)}\leq\,\rho U(y)\text{,}\\
		\label{38}	\lrr{z}{\Hess U(x)z}\leq\, c\norm{z}^2U(x)\text{.}
	\end{align}
	Then it holds for all $t\in\cc{0,T}$ that
	\begin{equation}
		\begin{aligned}
			\E\cc{U(Y^N_t)}\leq&\;\E\bigcc{U(Y_0^N)}\exp\biggoo{\Bigoo{\rho 		+2c^{p+3}\bigoo{\tfrac{T}{N}}^{\frac{15}{16}}		+32c^{3p+4}\bigoo{\tfrac{T}{N}}^{\frac{13}{32}}	}t}\\
			&+\int_0^t\bigoo{\tfrac{T}{N}}^{\frac{13}{32}}\Bigoo{32c^{3p+4}+\tfrac{4^{p+3}}{2}c^{p^2+4p+4}\bigoo{\tfrac{T}{N}}^p}\\
			&\hspace{5mm}\cdot\exp\biggoo{
				\Bigoo{\rho 		+2c^{p+3}\bigoo{\tfrac{T}{N}}^{\frac{15}{16}}		+32c^{3p+4}\bigoo{\tfrac{T}{N}}^{\frac{13}{32}}	}(t-s)}ds
			\text{.}
		\end{aligned}
	\end{equation}
\end{lemma}
\noindent
\textit{Proof of Lemma \ref{lemma 3}.} 
Let $\lrfloor{\cdot}_{N}\colon\cc{0,T}\to\cc{0,T}$ be the function which satisfies for all $t\in\oc{0,T}$ that $\lrfloor{t}_{N}=\sup(\gk{0,\frac{T}{N},\frac{2T}{N},\dots,T}\cap\co{0,t})$ and $\lrfloor{0}_{N}=0$. 
Let $C,\overline{C}\in\oo{0,\infty}$ satisfy that 
$C=\rho 		+2c^{p+3}\bigoo{\tfrac{T}{N}}^{\frac{15}{16}}		+32c^{3p+4}\bigoo{\tfrac{T}{N}}^{\frac{13}{32}}	$ and
$\overline{C}=\bigoo{\tfrac{T}{N}}^{\frac{13}{32}}\bigcc{32c^{3p+4}+\tfrac{4^{p+3}}{2}c^{p^2+4p+4}\bigoo{\tfrac{T}{N}}^p}$.
Taylor's theorem and (\ref{34}) imply for all $t\in\cc{0,T}$ that
\begin{equation}
	\begin{aligned}\label{50}
		&U\oo{Y_t^N}-U\oo{Y^N_{\lrfloor{t}_N}}\\
		=\;&\lrr{\nabla U(Y^N_{\lrfloor{t}_N})}{Y_t^N-Y^N_{\lrfloor{t}_N}}		
		+\tfrace \lrr{Y_t^N-Y^N_{\lrfloor{t}_N}}			{\Hess U(Y^N_{\lrfloor{t}_N})		(Y_t^N-Y^N_{\lrfloor{t}_N})}	\\
		&\hspace{5mm}+\tint{0}{1}\lrr{Y_t^N-Y^N_{\lrfloor{t}_N}}
		{\Bigoo{\Hess U\bigoo{\oo{1-\lambda}Y^N_{\lrfloor{t}_N}+\lambda Y^N_t}-\Hess U \oo{Y^N_{\lrfloor{t}_N}}}\oo{Y_t^N-Y^N_{\lrfloor{t}_N}}}\\
		&\hspace{10mm}\cdot\oo{1-\lambda}d\lambda\\
		\leq\;&\mathbbm{1}_{\norm{Y^N_{\lrfloor{t}_N}}\leq e^{\sqrt{\bb{\log(N/T)}}}}
		\cdot\biggoo{\lrr{\nabla U(Y^N_{\lrfloor{t}_N})	}
			{\mu(Y^N_{\lrfloor{t}_N})(t-{\lrfloor{t}_N})
				+\sigma(Y^N_{\lrfloor{t}_N})\Pi^N(W_t-W_{\lrfloor{t}_N})}\\
			&\hspace{30mm}+\tfrace\Biglr{\mu(Y^N_{\lrfloor{t}_N})(t-{\lrfloor{t}_N})
				+\sigma(Y^N_{\lrfloor{t}_N})\Pi^N(W_t-W_{\lrfloor{t}_N})}{\Hess U(Y^N_{\lrfloor{t}_N})\\		
				&\hspace{40mm}\cdot\Bigoo{\mu(Y^N_{\lrfloor{t}_N})(t-{\lrfloor{t}_N})
					+\sigma(Y^N_{\lrfloor{t}_N})\Pi^N(W_t-W_{\lrfloor{t}_N})}}
			\,}\\
		&\hspace{5mm}+\tint{0}{1}c\norm{Y_t^N-Y^N_{\lrfloor{t}_N}}^2
		\Bigoo{1+U(Y^N_{\lrfloor{t}_N})
			+\bignorm{Y_t^N-Y^N_{\lrfloor{t}_N}}^p}
		\cdot\bignorm{Y_t^N-Y^N_{\lrfloor{t}_N}}(1-\lambda) \,d\lambda\text{.}
	\end{aligned}
\end{equation}
The fact that $\forall j\neq k\in\gk{1,\dots,m}\colon\E\bigcc{\bigoo{\Pi^N(W_t-W_{\lrfloor{t}_N})}_j\bigoo{\Pi^N(W_t-W_{\lrfloor{t}_N})}_k|\F_{\lrfloor{t}_N}}=0$ a.s.
and $\forall k\in\gk{1,\dots,m}\colon \E\bigcc{\bigoo{\Pi^N\oo{W_t-W_{\lrfloor{t}_N}}}_k^2|\F_{\lrfloor{t}_N}}\leq\oo{t-\lrfloor{t}_N}$ a.s. show that for all $t\in\cc{0,T}$ it holds a.s. that
\begin{equation}
	\begin{aligned}
		\label{51}
		&\E\Bigcc{\lrr{\sigma(Y^N_{\lrfloor{t}_N})\Pi^N(W_t-W_{\lrfloor{t}_N})}
			{\Hess U(Y^N_{\lrfloor{t}_N})
				\sigma(Y^N_{\lrfloor{t}_N})\Pi^N(W_t-W_{\lrfloor{t}_N})}
			\Big|\F_{\lrfloor{t}_N}}\\
		=\;&\E\Biggcc{\sum_{i=1}^d\sum_{k=1}^m\sum_{j=1}^m
			\bigoo{\sigma(Y^N_{\lrfloor{t}_N})}_{ij}		\bigoo{\Pi^N(W_t-W_{\lrfloor{t}_N})}_j\\
			&\hspace{25mm}\cdot\bigoo{\Hess U(Y^N_{\lrfloor{t}_N})\sigma(Y^N_{\lrfloor{t}_N})}_{ik}		\bigoo{\Pi^N(W_t-W_{\lrfloor{t}_N})}_k
			\Bigg|\F_{\lrfloor{t}_N}}\\
		=\;&\sum_{i=1}^d\sum_{k=1}^m\sum_{j=1}^m
		\bigoo{\sigma(Y^N_{\lrfloor{t}_N})}_{ij}		\bigoo{\Hess U(Y^N_{\lrfloor{t}_N})\sigma(Y^N_{\lrfloor{t}_N})}_{ik}\\
		&\hspace{25mm}\cdot	\E\Bigcc{\bigoo{\Pi^N(W_t-W_{\lrfloor{t}_N})}_j\bigoo{\Pi^N(W_t-W_{\lrfloor{t}_N})}_k
			\Big|\F_{\lrfloor{t}_N}}\\
	=\;&\sum_{i=1}^d\sum_{k=1}^m
	\bigoo{\sigma(Y^N_{\lrfloor{t}_N})}_{ik}		\bigoo{\Hess U(Y^N_{\lrfloor{t}_N})\sigma(Y^N_{\lrfloor{t}_N})}_{ik}		\E\Bigcc{\bigoo{\Pi^N(W_t-W_{\lrfloor{t}_N})}_k^2
		\Big|\F_{\lrfloor{t}_N}}\\
	\leq\;&\lrrf{\sigma(Y^N_{\lrfloor{t}_N})}
	{\Hess U(Y^N_{\lrfloor{t}_N})
		\bigoo{\sigma(Y^N_{\lrfloor{t}_N})}}
	(t-\lrfloor{t}_N)\text{.}
\end{aligned}
\end{equation}
Furthermore, for all $t\in\cc{0,T}$ it holds a.s. that
\begin{equation}
\begin{aligned}
	\E\,\biggcc{&\Biglr{\nabla U(Y^N_{\lrfloor{t}_N})	}{\sigma(Y^N_{\lrfloor{t}_N})\Pi^N(W_t-W_{\lrfloor{t}_N})}\\
		+\tfrace&\Biglr{\mu(Y^N_{\lrfloor{t}_N})(t-{\lrfloor{t}_N})}
		{\Hess U(Y^N_{\lrfloor{t}_N})\sigma(Y^N_{\lrfloor{t}_N})\Pi^N(W_t-W_{\lrfloor{t}_N})}\\
		+\tfrace&\Biglr{\sigma(Y^N_{\lrfloor{t}_N})\Pi^N(W_t-W_{\lrfloor{t}_N})}
		{\Hess U(Y^N_{\lrfloor{t}_N})\mu(Y^N_{\lrfloor{t}_N})(t-{\lrfloor{t}_N})}
		\bigg|\mathcal{F}_{\lrfloor{t}_N}}\\
	=\hspace{7mm}&\Biglr{\nabla U(Y^N_{\lrfloor{t}_N})}
	{\sigma(Y^N_{\lrfloor{t}_N})\E\Bigcc{\Pi^N(W_t-W_{\lrfloor{t}_N})\Big|\mathcal{F}_{\lrfloor{t}_N}}}\\
	+\tfrace&\Biglr{\mu(Y^N_{\lrfloor{t}_N})}
	{\Hess U(Y^N_{\lrfloor{t}_N})\sigma(Y^N_{\lrfloor{t}_N})\E\Bigcc{\Pi^N(W_t-W_{\lrfloor{t}_N})\Big|\mathcal{F}_{\lrfloor{t}_N}}}
	(t-{\lrfloor{t}_N})\\
	+\tfrace&\Biglr{\sigma(Y^N_{\lrfloor{t}_N})\E\Bigcc{\Pi^N(W_t-W_{\lrfloor{t}_N})\Big|\mathcal{F}_{\lrfloor{t}_N}}}
	{\Hess U(Y^N_{\lrfloor{t}_N})\mu(Y^N_{\lrfloor{t}_N})}(t-{\lrfloor{t}_N})
	=0\hspace{0mm}\text{.}
\end{aligned}
\end{equation}
This, (\ref{50}), and (\ref{51}) show that for all $t\in\cc{0,T}$ it holds a.s. that
\begin{equation}
\begin{aligned}
	\label{52}
	&\E\lrcc{U\oo{Y_t^N}
		-U\oo{Y^N_{\lrfloor{t}_N}}|\F_{\lrfloor{t}_N}}\\
	\leq\;&\E\biggcc{
		\mathbbm{1}_{\norm{Y^N_{\lrfloor{t}_N}}\leq e^{\sqrt{\bb{\log(N/T)}}}}
		\cdot\Bigoo{\lrr{\nabla U(Y^N_{\lrfloor{t}_N})}		{\mu(Y^N_{\lrfloor{t}_N})}\\
		&\hspace{10mm}+\tfrace\Biglr{\sigma(Y^N_{\lrfloor{t}_N})}{\Hess U(Y^N_{\lrfloor{t}_N})	\sigma(Y^N_{\lrfloor{t}_N})}_F
		}			(t-{\lrfloor{t}_N})\\
		&\hspace{5mm}+	\mathbbm{1}_{\norm{Y^N_{\lrfloor{t}_N}}\leq e^{\sqrt{\bb{\log(N/T)}}}}\cdot\tfrace\biglr{\mu(Y^N_{\lrfloor{t}_N})}			{\Hess U(Y^N_{\lrfloor{t}_N})	
			\mu(Y^N_{\lrfloor{t}_N})}				\,(t-{\lrfloor{t}_N})^2\\
		&\hspace{5mm}+\tint{0}{1}c\bigoo{1+U(Y^N_{\lrfloor{t}_N})}\norm{Y_t-Y_{\lrfloor{t}_N}}^3(1-\lambda)\,d\lambda
		+\tint{0}{1}c\norm{Y_t-Y_{\lrfloor{t}_N}}^{p+3}(1-\lambda)\,d\lambda
		\bigg|\F_{\lrfloor{t}_N}}\text{.}
\end{aligned}
\end{equation}
Furthermore, (\ref{112}) and (\ref{113}) show for all $t\in\cc{0,T}$ that
\begin{equation}
\label{54}
\maxofv{\bignorm{\mu(Y^N_{\lrfloor{t}_N})}}{\bignormf{\sigma(Y^N_{\lrfloor{t}_N})}}
\,\leq\, c\bigoo{1+\norm{Y^N_{\lrfloor{t}_N}}^p}
\,\leq\, 2ce^{p\sqrt{\bb{\log(N/T)}}}
\,\leq\, 2c^{p+1}\bigoo{\tfrac{N}{T}}^{\frac{1}{32}}\text{.}
\end{equation}
This, the triangle inequality, (\ref{42}), and the fact that $T/N\in\oc{0,1}$ yield that for all $t\in\cc{0,T}$, $q\in\gk{3,p+3}$ it holds a.s. that 
\begin{equation}
\begin{aligned}
	&\E\Bigcc{\bignorm{Y_t^N-Y^N_{\lrfloor{t}_N}}^q\Big|\F_{\lrfloor{t}_N}}\\
	=\;&\E\Bigcc{\mathbbm{1}_{\norm{Y^N_{\lrfloor{t}_N}}\leq e^{\sqrt{\bb{\log(N/T)}}}}\cdot\bignorm{\mu(Y^N_{\lrfloor{t}_N})(t-{\lrfloor{t}_N})
			+\sigma(Y^N_{\lrfloor{t}_N})\Pi^N(W_t-W_{\lrfloor{t}_N})}^q\Big|\F_{\lrfloor{t}_N}}\\
	\leq\;&\E\Bigcc{\mathbbm{1}_{\norm{Y^N_{\lrfloor{t}_N}}\leq e^{\sqrt{\bb{\log(N/T)}}}}\cdot\Bigoo{\bignorm{\mu(Y^N_{\lrfloor{t}_N})}(t-{\lrfloor{t}_N})\\
		&\hspace{5mm}+\bignorm{\sigma(Y^N_{\lrfloor{t}_N})}_{L(\Rd,\Rm)}\,\bignorm{\Pi^N(W_t-W_{\lrfloor{t}_N})}}^q\Big|\F_{\lrfloor{t}_N}}\\
	\leq\;&\E\Bigcc{\mathbbm{1}_{\norm{Y^N_{\lrfloor{t}_N}}\leq e^{\sqrt{\bb{\log(N/T)}}}}
		\cdot\Bigoo{\maxofv{\bignorm{\mu(Y^N_{\lrfloor{t}_N})}\!\!}{\!\!\bignorm{\sigma(Y^N_{\lrfloor{t}_N})}_{L(\Rd,\Rm)}}}^q\\
		&\hspace{5mm}\cdot\Bigoo{(t-{\lrfloor{t}_N})
			+\bignorm{\Pi^N(W_t-W_{\lrfloor{t}_N})}}^q\Big|\F_{\lrfloor{t}_N}}\\
	\leq\;&2^{q}c^{q(p+1)}\bigoo{\tfrac{N}{T}}^{\frac{q}{32}}
	\Bigcc{\oo{t-\lrfloor{t}_N}+\Bigoo{\E\Bigcc{\bignorm{W_t-W_{\lrfloor{t}_N}}^q\Big|\F_{\lrfloor{t}_N}}}^{\frac{1}{q}}}^q\\
	\leq\;&2^{q}c^{q(p+1)}\bigoo{\tfrac{N}{T}}^{\frac{q}{32}}
	\Bigoo{\oo{t-\lrfloor{t}_N}^{\frac{1}{2}}+q\oo{t-\lrfloor{t}_N}^{\frac{1}{2}}}^q\\
	=	\;&2^{q}(q+1)c^{q(p+1)}\bigoo{\tfrac{N}{T}}^{\frac{q}{32}}
	\oo{t-\lrfloor{t}_N}^{\frac{q}{2}-1}\oo{t-\lrfloor{t}_N}\\
	\leq\;&4^{q}c^{q(p+1)}
	\bigoo{\tfrac{T}{N}}^{\frac{15}{32}q-1}\oo{t-\lrfloor{t}_N}\text{.}
\end{aligned}
\end{equation}
This, (\ref{33}), (\ref{38}), (\ref{52}), (\ref{54}), the fact that $\int_0^1(1-\lambda)\,d\lambda=\frace$, and $T/N\in\oo{0,1}$ imply that for all $t\in\cc{0,T}$ it holds a.s. that
\begin{equation}
\begin{aligned}
	\label{46}
	&\E\lrcc{U\oo{Y_t^N}
		-U\oo{Y^N_{\lrfloor{t}_N}}|\F_{\lrfloor{t}_N}}\\
	\leq\;&\mathbbm{1}_{\norm{Y^N_{\lrfloor{t}_N}}\leq e^{\sqrt{\bb{\log(N/T)}}}}\!\cdot\!\Bigoo{
		\rho U(Y^N_{\lrfloor{t}_N})	\oo{t-\lrfloor{t}_N}+\tfrace c\bignorm{\mu(Y^N_{\lrfloor{t}_N})}^2U(Y^N_{\lrfloor{t}_N})(t-{\lrfloor{t}_N})^2}\\
	&\hspace{5mm}+32c^{3p+4}\bigoo{\tfrac{T}{N}}^{\frac{13}{32}}\bigoo{1+U(Y^N_{\lrfloor{t}_N})}(t-\lrfloor{t}_N)
	+\tfrac{4^{p+3}}{2}c^{p^2+4p+4}\bigoo{\tfrac{T}{N}}^{p+\frac{13}{32}}\oo{t-\lrfloor{t}_N}\\
	\leq\;&\Bigoo{\rho 		+2c^{p+3}\bigoo{\tfrac{T}{N}}^{\frac{15}{16}}		+32c^{3p+4}\bigoo{\tfrac{T}{N}}^{\frac{13}{32}}	}		
	U(Y^N_{\lrfloor{t}_N})	\oo{t-\lrfloor{t}_N}\\
	&\hspace{5mm}+\bigoo{\tfrac{T}{N}}^{\frac{13}{32}}\Bigoo{32c^{3p+4}+\tfrac{4^{p+3}}{2}c^{p^2+4p+4}\bigoo{\tfrac{T}{N}}^p}\oo{t-\lrfloor{t}_N}\\
	=\;&CU(Y^N_{\lrfloor{t}_N})	\oo{t-\lrfloor{t}_N}+\overline{C}\oo{t-\lrfloor{t}_N}\text{.}
\end{aligned}
\end{equation}
Furthermore (\ref{46}) shows for all $t\in\cc{0,T}$ that
\begin{equation}
\begin{aligned}
	&\E\bigcc{U(Y^N_t)}\\
	\leq\;&\E\bigcc{U(Y^N_{\lrfloor{t}_N})}\bigoo{1+C(t-\lrfloor{t}_N)}		+\overline{C}(t-\lrfloor{t}_N)\\
	\leq\;&\E\bigcc{U(Y^N_{\lrfloor{t}_N})}\exp\bigoo{C(t-\lrfloor{t}_N)}
	+\int_{\lrfloor{t}_N}^{t}\overline{C}	\exp\bigoo{C(t-s)}ds\\
	\leq\;&\E\bigcc{U(Y_{\lrfloor{\lrfloor{t}_N}_N}^N)}	\exp\bigoo{C\oo{t-\lrfloor{\lrfloor{t}_N}_N}}\\
	&+\int_{\lrfloor{\lrfloor{t}_N}_N}^{\lrfloor{t}_N}\overline{C}\exp\bigoo{C(t-s)}
	+\int_{\lrfloor{t}_N}^{t}\overline{C}	\exp\bigoo{C(t-s)}ds\\
	=\;&\E\bigcc{U(Y_{\lrfloor{\lrfloor{t}_N}_N}^N)}	\exp\bigoo{C\oo{t-\lrfloor{\lrfloor{t}_N}_N}}
	+\int_{\lrfloor{\lrfloor{t}_N}_N}^{t}\overline{C}	\exp\bigoo{C(t-s)}ds\\
	\vdots\;&\\
	\leq\;&\E\bigcc{U(Y_0^N)}	\exp\bigoo{Ct}	+\int_0^{t}\overline{C}	\exp\bigoo{C(t-s)}ds\text{.}
\end{aligned}
\end{equation}
The proof of Lemma \ref{lemma 3} is thus completed. $ \hfill\square $\\
\section{Strong convergence rate}
The following theorem is our main result 
and states that stopped Brownian-increment tamed Euler approximations
converge with strong rate $1/2$ to the exact solution.
Theorem \ref{lemma 4} is motivated by \cite{hutzenthaler2020perturbation} which proves the
respective result for increment-tamed Euler approximations.
\begin{theorem}
	\label{lemma 4}
	Let $d,m\in\N$, 
	$T,r\in\iooi$, 
	$c\in\co{T^{\frac{1}{32}},\infty}$,
	$q_0,q_1\in\ioci$, 
	$\rho\in\icoi$, 
	$p,r\in\co{2,\infty}$ with $\frac{1}{p}+\frac{1}{q_0}+\frac{1}{q_1}=\frac{1}{r}$,
	let $\mu\colon\Rd\to\Rd$, 
	$\sigma\colon\Rd\to\Rdm$
	be locally Lipschitz continuous functions with at most polynomially growing Lipschitz constants,
	let $U\in C^2(\Rd,\co{0,\infty})$, 
	let $\overline{U}\colon\Rd\to\R$ be measurable function, 
	assume that $\Hess (U)$, $\overline{U}$ are at most polynomially growing functions,
	for all $N\in\N$ let $\Pi^N\in C^2(\Rm,\Rm)$ be the function which satisfies for all $x=\oo{x_1,\dots,x_m}\in\Rm$ that
	\begin{equation}
		\Pi^N(x)=\lroo{x_i\,\exp\lroo{-\frac{Nx_i^4}{T}}}_{i\in\gk{1,\dots,m}}\text{,}
	\end{equation}
	let $(\Omega,\F,\P,(\mathbb{F}_t)_{t\in\cc{0,T}})$ be a filtered probability space, 
	let $W\colon\icc{0}{T}\times\Omega\to\R^m$ be a standard $(\mathbb{F}_t)_{t\in\cc{0,T}}$-Brownian motion with continuous sample paths,
	for all $N\in\N$ let $X\colon\icc{0}{T}\times\Omega\to\Rd$ and $Y^N\colon\cc{0,T}\times\Omega\to\Rd$ be 
	adapted stochastic processes with continuous sample paths which satisfy for all $k\in\gk{0,1,\dots,N-1}$, $t\in\bigcc{\frac{kT}{N},\frac{(k+1)T}{N}}$ that
	\begin{align}
		X_t&=X_0+\int_0^t\mu(X_s)ds+\int_0^t\sigma(X_s)dW_s\\
		Y_t^N&=Y^N_{\frac{kT}{N}}+\mathbbm{1}_{\norm{Y^N_{\frac{kT}{N}}}\leq e^{\sqrt{\bb{\log(T/N)}}}}
		\cdot\Bigoo{\mu\bigoo{Y^N_{\frac{kT}{N}}}\bigoo{t-\tfrac{kT}{N}}
			+\sigma\bigoo{Y^N_{\frac{kT}{N}}}\Pi^N\bigoo{W_t-W_{\frac{kT}{N}}}}\label{eq:YN}
	\end{align}
	and $Y^N_0=X_0$, assume $\E\bigcc{e^{U(X_0)}}<\infty$
	and assume for all $t\in\cc{0,T}$, $N\in\N$, $k\in\gk{0,1,\dots,N-1}$, $x,y\in\Rd$ that
	\begin{align}
		\label{48}	\frac{\lr{x-y}{\mu(x)-\mu(y)}+\frac{(p-1)(1+1/c)}{2}\normf{\sigma(x)-\sigma(y)}^2}{\norm{x-y}^2}
		&\leq\!c\!+\!\frac{\bb{U(x)}\+ \bb{U(y)}}{2q_0Te^{\rho T}}\+ \frac{\bb{\overline{U}(x)}\+ \bb{\overline{U}(y)}}{2q_1e^{\rho T}}\text{,}\\
		\biglr{(\nabla U)(x)}{\mu(x)}+\tfrace\biglr{\sigma(x)}{(\Hess U)(x)\sigma(x)}_{\tF}
		+\tfrace&\norm{\oo{\sigma(x)}^*(\nabla U)(x)}^2+\overline{U}(x)\leq\rho\cdot U(x)\\
		\label{0152}	\tfrac{1}{c}\norm{x}^{1/c}&\leq1+\bb{U(x)}\text{.}
	\end{align}
	Then there exists a real number $C\in\icoi$ such that for all $N\in\N$ it holds that
	\begin{equation}
		\label{71}
		\sup_{k\in\gk{0,1,\dots,N}}\lrnorm{X_{\frac{kT}{N}}-Y^N_{\frac{kT}{N}}}_{L^r(\Omega;\Rd)}\leq C\bigoo{\tfrac{T}{N}}^{1/2}\text{.}
	\end{equation}
\end{theorem}
\noindent
\textit{Proof of Theorem \ref{lemma 4}.} 
Let $q\in\oc{0,\infty}$ satisfy that $\frac{1}{q}=\frac{1}{q_0}+\frac{1}{q_1}$.
For all $N\in\N$ let $\tau^N\colon\Omega\to\cc{0,T}$ be the function which satisfies that
$\tau^N=\inf\bigoo{\biggk{\frac{kT}{N}\in\cc{0,T}\colon k\in\gk{0,1,\dots,N},\, \norm{Y_{\frac{kT}{N}}^N}>e^{\sqrt{\bb{\log(N/T)}}}}\cup\gk{T}}$.
Lemma \ref{lemma 2} and [Corollary 2.4,\cite{cox2022local}] imply that there exists $C_1\in\R$ such that it holds for all $N\in\N$ that
\begin{equation}
	\label{49}
	\begin{aligned}
		&\lrcc{\sup_{t\in\cc{0,T}}\E\lrcc{\exp\biggoo{e^{-\rho t}\bb{U\oo{X_t}}
					+{\textstyle\int\limits_0^t}e^{-\rho u}\bb{\overline{U}\oo{X_u}}\,du}}}\\
		\cdot&\lrcc{\sup_{t\in\cc{0,T}}\E\biggcc{\exp\biggoo{e^{-\rho (\minofv{t}{\tau^N})}\bb{U\oo{Y^N_t}}
					+{\textstyle\int\limits_0^{\minofv{t}{\tau^N}}} e^{-\rho u}\bb{\overline{U}\oo{Y^N_u}}\,du}}}\leq C_1\text{.}
	\end{aligned}
\end{equation}
This, (\ref{0152}), and the fact that $\forall x\in\Rd,t\in\cc{0,T}\colon U(x)\leq e^{\rho T}\exp(-e^{-\rho t}U(x))$ 
and $\norm{x}^{2r}\leq(c[1+U(x)])^{2rc}\leq\sup_{z\in\R}\frac{(c(1+e^{\rho T}z))^{2rc}}{\exp(z)}\exp(e^{-\rho t}U(x))$ imply that there exists $C_2\in\R$ such that it holds for all $N\in\N$ that
\begin{equation}
	\label{151}
	\sup_{t\in\cc{0,T}}\bignorm{U(Y^N_t)}_{L\oo{\Omega;\R}}^{\frac{1}{2p}}
	+\sup_{t\in\cc{0,T}}\Bigcc{\bignorm{X_t}_{L^{2r}\oo{\Omega;\Rd}}+\bignorm{Y^N_t}_{L^{2r}\oo{\Omega;\Rd}}}\leq C_2\text{.}
\end{equation}
Combining Corollary 2.12  in \cite{hutzenthaler2020perturbation} (applied for all $N\in\N$, $t\in\cc{0,T}$ with $\epsilon\curvearrowleft \frac{1}{c}$, $\tau\curvearrowleft\minofv{t\!}{\!\tau^N}$, $Y\curvearrowleft Y^N$, $a\curvearrowleft \bigoo{\mu(Y^N_{\lrfloor{s}_N})\+ \tfrace\sigma(Y^N_{\lrfloor{s}_N})\Delta\Pi^N\oo{W_s\!-\!W_{\lrfloor{s}_N}}}_{s\in\cc{0,T}}$, $b\curvearrowleft \bigoo{\sigma(Y^N_{\lrfloor{s}_N})D\Pi^N\oo{W_s\- W_{\lrfloor{s}_N}}}_{s\in\cc{0,T}}$, $\delta\curvearrowleft 1$, $\rho\curvearrowleft 1$) and Lemma \ref{lemma 5} (applied for every $N\in\N$, $t\in\cc{0,T}$) shows for all $N\in\N$ that
\begin{equation}
	\label{74}
	\begin{aligned}
		&\sup_{t\in\cc{0,T}}\bignorm{X_{\minofv{t}{\tau^N}}-Y^N_{\minofv{t}{\tau^N}}}_{L^r(\Omega;\Rd)}\\
		\leq\,&\lrnorm{\exp\!\lroo{{\textstyle\int\limits_0^{\tau^N}}\!
				\lrcc{\frac{\lr{X_s\!-\!Y_s^N\!}{\!\mu(X_s)-\mu(Y_s^N)}
						+\frac{(p-1)(1+1/c)}{2}\normf{\sigma(X_s)\!-\!\sigma(Y_s^N)}^2}
					{\norm{X_s-Y_s^N}^2}+\tfrac{3}{2}-\tfrac{2}{p}}^+\!\!ds}}_{L^q(\Omega;\R)}\\
		&\cdot\biggoo{\Bigoo{\tint{0}{T}\bignorm{\mu(Y^N_{\lrfloor{s}_N})\+ \tfrace\sigma(Y^N_{\lrfloor{s}_N})\Delta\Pi^N\oo{W_s\!-\!W_{\lrfloor{s}_N}}-\mu(Y_s^N)}^pds}^{\frac{1}{p}}\\
			&\hspace{5mm}+\sqrt{(p\- 1)(1\+ c)}\cdot\Bigoo{\tint{0}{T}\bignormf{\sigma(Y^N_{\lrfloor{s}_N})D\Pi^N\oo{W_s\!-\!W_{\lrfloor{s}_N}}\- \sigma(Y_s^N)}^pds}^{\frac{1}{p}}}\text{.}
	\end{aligned}
\end{equation}
Assumption (\ref{48}), Hölder's inequality, and Jensen's inequality imply for all $N\in\N$ that
\begin{equation}
	\begin{aligned}
		&\lrnorm{\exp\!\lroo{{\textstyle\int\limits_0^{\tau^N}}\!
				\lrcc{\frac{\lr{X_s\!-\!Y_s^N\!}{\!\mu(X_s)-\mu(Y_s^N)}
						+\frac{(p-1)(1+1/c)}{2}\normf{\sigma(X_s)\!-\!\sigma(Y_s^N)}^2}
					{\norm{X_s-Y_s^N}^2}}^+\!\!ds}}_{L^q(\Omega;\R)}e^{-cT}\\
		&\leq\lrnorm{\exp\!\lroo{{\textstyle\int\limits_0^{\tau^N}} \frac{\bb{U(X_s)}+\bb{U(Y_s^N)}}{2q_0Te^{\rho T}}
				+\frac{\bb{\overline{U}(X_s)}+\bb{\overline{U}(Y_s^N)}}{2q_1e^{\rho T}}		ds}}_{L^q(\Omega;\R)}\\
		&\leq\!\sup_{s\in\cc{0,T}}\!\biggnorm{\frac{\bb{U(X_s)}}{2q_0e^{\rho T}}}_{L^{2q_0}(\Omega;\R)}
		\biggnorm{\textstyle\int\limits_0^{\tau^N}\!\displaystyle\frac{\bb{\overline{U}(X_u)}}{2q_1e^{\rho T}}du}_{L^{2q_1}(\Omega;\R)}
		\biggnorm{\frac{\bb{U(Y^N_s)}}{2q_0e^{\rho T}}}_{L^{2q_0}(\Omega;\R)}
		\biggnorm{\textstyle\int\limits_0^{\tau^N}\displaystyle\!\frac{\bb{\overline{U}(Y^N_u)}}{2q_1e^{\rho T}}du}_{L^{2q_1}(\Omega;\R)}
		\text{.}
	\end{aligned}
\end{equation}
This, \ref{023}. in Lemma \ref{lemma 6} (applied for every $N\in\N$, $t\in\cc{0,T}$), (\ref{49}), (\ref{74}), Lemma \ref{lemma 1} (applied for every $N\in\N$, $s\in\cc{0,T}$ with $h\curvearrowleft T/N$, $W\curvearrowleft (W_{\lrfloor{s}_N+t}-W_{\lrfloor{s}_N})_{t\in\cc{0,T/N}}$), and the triangle inequality imply that there exists $C_3\in\R$ such that it holds for all $N\in\N$ that
\begin{equation}
	\begin{aligned}
		&\sup_{t\in\cc{0,T}}\bignorm{X_{\minofv{t}{\tau^N}}-Y^N_{\minofv{t}{\tau^N}}}_{L^r\oo{\Omega;\Rd}}\\
		\leq&\,\lrnorm{\exp\!\lroo{{\textstyle\int\limits_0^{\tau^N}}\!
				\lrcc{\frac{\lr{X_s\!-\!Y_s^N\!}{\!\mu(X_s)-\mu(Y_s^N)}
						+\frac{(p-1)(1+1/c)}{2}\normf{\sigma(X_s)\!-\!\sigma(Y_s^N)}^2}
					{\norm{X_s-Y_s^N}^2}+\tfrac{3}{2}-\tfrac{2}{p}}^+\!\!ds}}_{L^q(\Omega;\R)}\\
		&\cdot\biggoo{\Bigoo{\tint{0}{T}\bignorm{\mu(Y^N_{\lrfloor{s}_N})\+ \tfrace\sigma(Y^N_{\lrfloor{s}_N})\Delta\Pi^N\oo{W_s\!-\!W_{\lrfloor{s}_N}}-\mu(Y_s^N)}^pds}^{\frac{1}{p}}\\
			&\hspace{5mm}+\sqrt{(p\- 1)(1\+ c)}\cdot\Bigoo{\tint{0}{T}\bignormf{\sigma(Y^N_{\lrfloor{s}_N})D\Pi^N\oo{W_s\!-\!W_{\lrfloor{s}_N}}\- \sigma(Y_s^N)}^pds}^{\frac{1}{p}}}\\
		\leq&\,\sup_{s\in\cc{0,T}}\!\Biggoo{\biggnorm{\frac{\bb{U(X_s)}}{2q_0e^{\rho T}}}_{L^{2q_0}(\Omega;\R)}
			\biggnorm{\textstyle\int\limits_0^{\tau^N}\!\displaystyle\frac{\bb{\overline{U}(X_u)}}{2q_1e^{\rho T}}du}_{L^{2q_1}(\Omega;\R)}\\
			&\hspace{10mm}\cdot\biggnorm{\frac{\bb{U(Y^N_s)}}{2q_0e^{\rho T}}}_{L^{2q_0}(\Omega;\R)}
			\biggnorm{\textstyle\int\limits_0^{\tau^N}\displaystyle\!\frac{\bb{\overline{U}(Y^N_u)}}{2q_1e^{\rho T}}du}_{L^{2q_1}(\Omega;\R)}}\\
		&\cdot\Bigoo{\bignorm{\mu(Y^N_{\lrfloor{t}_N})\- \mu(Y_t^N)}_{L^p(\Omega;\Rd)}
			\+ \tfrace\bignorm{\sigma(Y^N_{\lrfloor{t}_N})\Delta\Pi^N\oo{W_t\!-\!W_{\lrfloor{t}_N}}}_{L^p(\Omega;\Rd)}+\sqrt{(p\- 1)(1\+ c)}\\
			&\hspace{5mm}\cdot\bignorm{\sigma(Y^N_{\lrfloor{t}_N})\bigoo{D\Pi^N\oo{W_t\!-\!W_{\lrfloor{t}_N}}-I_{\Rdm}}}_{L^p(\Omega;\Rdm)}+ \bignorm{\sigma(Y^N_{\lrfloor{t}_N})\-\sigma(Y_t^N)}_{L^p(\Omega;\Rdm)}}\\
		\leq&\,\sup_{s\in\cc{0,T}}\!\biggnorm{\frac{\bb{U(X_s)}}{2q_0e^{\rho (\minofv{s\!}{\!\tau^N})}}+
			\textstyle\int\limits_0^s\!\displaystyle\frac{\bb{\overline{U}(X_u)}}{2q_0e^{\rho u}}du}_{L^{2q_0}(\Omega;\R)}
		\cdot\biggnorm{\frac{\bb{U(X_s)}}{2q_1e^{\rho (\minofv{s\!}{\!\tau^N})}}+
			\textstyle\int\limits_0^s\!\displaystyle\frac{\bb{\overline{U}(X_u)}}{2q_1e^{\rho u}}du}_{L^{2q_1}(\Omega;\R)}\\
		&\cdot\biggnorm{\frac{\bb{U(Y^N_s)}}{2q_0e^{\rho T}}+
			\textstyle\int\limits_0^{\minofv{s\!}{\!\tau^N}}\displaystyle\!\frac{\bb{\overline{U}(Y^N_u)}}{2q_0e^{\rho T}}du}_{L^{2q_0}(\Omega;\R)}
		\cdot\biggnorm{\frac{\bb{U(Y^N_s)}}{2q_1e^{\rho T}}+
			\textstyle\int\limits_0^{\minofv{s\!}{\!\tau^N}}\displaystyle\!\frac{\bb{\overline{U}(Y^N_u)}}{2q_1e^{\rho T}}du}_{L^{2q_1}(\Omega;\R)}\\
		&\cdot\Bigoo{4c^{2p+2}3^{2p}\bigoo{\tfrac{N}{T}}^{\frac{1}{16}}	\Bigcc{2c^{p+1}\bigoo{\tfrac{T}{N}}^{\frac{7}{32}}
				\bigoo{T^{\frac{3}{4}}\+ \sqrt{m}}\+ 16\bigoo{\tfrac{T}{N}}^{\frace}\sqrt{m}}\\
			&\hspace{5mm}\+ \sqrt{(p\- 1)(1\+ c)}\cdot\Bigcc{52\tfrac{T}{N}\sqrt{m}\cdot 2c^{p+1}\bigoo{\tfrac{N}{T}}^{\frac{1}{32}}
				\+ 2c^{2p+2}3^{p}\bigoo{\tfrac{T}{N}}^{\frac{3}{16}}\bigoo{T^{\frac{3}{4}}\+ \sqrt{m}}}}\\
		\leq&\, C_3\bigoo{\tfrac{T}{N}}^{1/2}\text{.}
	\end{aligned}
\end{equation}
This, the triangle inequality, Hölder's inequality, and (\ref{151}) show for all $N\in\N$ that 
\begin{equation}
	\begin{aligned}
		\label{72}
		&\sup_{t\in\cc{0,T}}\bignorm{X_t-Y^N_t}_{L^r\oo{\Omega;\Rd}}\\
		=\,&\sup_{t\in\cc{0,T}}\Bignorm{\mathbbm{1}_{\gk{\tau^N<T}}\cdot\bigoo{X_t-Y^N_t}
			+\mathbbm{1}_{\gk{\tau^N\geq T}}\cdot\bigoo{X_{\minofv{t}{\tau^N}}-Y^N_{\minofv{t}{\tau^N}}}}_{L^r\oo{\Omega;\Rd}}\\
		\leq\,&\norm{\mathbbm{1}_{\tau^N<T}}_{L^{2r}\oo{\Omega;\R}}
		\sup_{t\in\cc{0,T}}\Bigcc{\bignorm{X_t}_{L^{2r}\oo{\Omega;\Rd}}+\bignorm{Y^N_t}_{L^{2r}\oo{\Omega;\Rd}}}
		\!+\!\sup_{t\in\cc{0,T}}\bignorm{X_{\minofv{t}{\tau^N}}-Y^N_{\minofv{t}{\tau^N}}}_{L^r\oo{\Omega;\Rd}}\\
		\leq\,& C_2\bigbb{\P\cc{\tau^N<T}}^{\frac{1}{2r}}
		+C_3\bigoo{\tfrac{T}{N}}^{1/2}\text{.}
	\end{aligned}
\end{equation}
Assumption, Markov's inequality, the fact that for all $x\in\icoi$ it holds that $\frac{1}{4!}x^4\leq e^x$, and (\ref{49}) imply for all $N\in\N$ that
\begin{equation}
	\begin{aligned}
		\label{73}
		&\P\cc{\tau^N<T}\\
		\leq\,&\P\Bigcc{\norm{Y^N_T}\geq \exp\bigoo{\sqrt{\bb{\log\oo{T/N}}}}}\\
		\leq\,&\P\biggcc{\frac{1+\bb{U(Y_T^N)}}{e^{\rho \tau^N}}
			\geq\frac{1}{c\,e^{\rho T}}\,\exp\Bigoo{\tfrac{\sqrt{\bb{\log\oo{T/N}}}}{c}}}\\
		\leq\,&\E\biggcc{\exp\biggoo{\frac{1+\bb{U(Y_T^N)}}{e^{\rho \tau^N}}}}
		\exp\biggoo{\!\!-\tfrac{1}{c\,e^{\rho T}}\exp\Bigoo{\tfrac{\sqrt{\bb{\log\oo{T/N}}}}{c}}}\\
		\leq\,&C_1 \exp\biggoo{1-\frac{\cc{\log\oo{T/N}}^2}{24c^5\,e^{\rho T}}}\text{.}
	\end{aligned}
\end{equation}
There exists $n_0\in\N$ such that for all $z\in\co{n_0,\infty}$ it holds that $\exp\oo{-\cc{\log(T/z)}^2/\oo{48rc^5\,e^{\rho T}}}\leq\sqrt{T/z}$.
This and (\ref{73}) show for all $N\in\N\cap\co{n_0,\infty}$ that
\begin{equation}
	\begin{aligned}
		&\bigbb{\P\cc{\tau^N<T}}^{\frac{1}{2r}}\\
		\leq\,&\biggbb{C_1 \exp\biggoo{1-\frac{\cc{\log\oo{T/N}}^2}{24c^5\,e^{\rho T}}}}^{\frac{1}{2r}}\\
		\leq\,&\Bigbb{C_1\exp\bigoo{1}}^{\frac{1}{2r}}
		\cdot\exp\biggoo{-\frac{\cc{\log\oo{T/N}}^2}{48rc^5\,e^{\rho T}}}\\
		\leq\,&\Bigbb{C_1\exp\bigoo{1}}^{\frac{1}{2r}}\bigoo{\tfrac{T}{N}}^{1/2}\text{.}
	\end{aligned}
\end{equation}
This, (\ref{151}), and (\ref{72}) imply (\ref{71}).
The proof of Lemma \ref{lemma 4} is thus completed. $ \hfill\square $

\subsubsection*{Acknowledgements}
This work has been funded by the Deutsche Forschungsgemeinschaft (DFG, German Research Foundation) via
RTG 2131 \textit{High-dimensional Phenomena in Probability -- Fluctuations and Discontinuity.}
%


\end{document}